\documentclass[12pt]{amsart}

\usepackage{amssymb,latexsym,amsmath,epsfig}
\usepackage[hyphens]{url} 
\usepackage[colorlinks,linkcolor=blue]{hyperref}
\usepackage[utf8]{inputenc}
\usepackage[T1]{fontenc}
\usepackage{graphicx}
\usepackage{mathrsfs}

 \hypersetup{colorlinks = true,linkcolor = blue,anchorcolor =red,citecolor = blue,filecolor = red,urlcolor = red,pdfauthor=author}

\def\pf{\begin{proof}}

	\usepackage{amssymb,textcomp, url}
	
\usepackage{geometry}

 \newtheorem{thm}{Theorem}[section]
 \newtheorem{prop}[thm]{Proposition}
 \newtheorem{lem}[thm]{Lemma}
 \newtheorem{cor}[thm]{Corollary}
\theoremstyle{definition}
 
 \newtheorem{defn}[thm]{Definition}

 \newtheorem{rem}[thm]{Remark}
 \numberwithin{equation}{section}

 \usepackage{t1enc,mathrsfs,amssymb,amsfonts,latexsym,pifont,fancybox}
	\newcommand{\bprop} {\begin{proposition}}
		\newcommand{\eprop} {\end{proposition}}
	\newcommand{\btheo} {\begin{theorem}}
		\newcommand{\etheo} {\end{theorem}}
	\newcommand{\blem} {\begin{lemma}}
		\newcommand{\elem} {\end{lemma}}
	\newcommand{\bcor} {\begin{corollary}}
		\newcommand{\ecor} {\end{corollary}}
	
	\newcounter{rea}
	\newcounter{red}
	\newcommand{\Be}{\begin{equation}}
	\newcommand{\Ee}{\end{equation}}
	\newcommand{\Bea}{\begin{eqnarray}}
	\newcommand{\Eea}{\end{eqnarray}}
	\newcommand{\Bes}{\begin{equation*}}
	\newcommand{\Ees}{\end{equation*}}
	\newcommand{\Beas}{\begin{eqnarray*}}
		\newcommand{\Eeas}{\end{eqnarray*}}
	\newcommand{\Ba}{\begin{array}}
		\newcommand{\Ea}{\end{array}}
	\def\R{\mathbb{R}}

	\scrollmode

	\title[Weighted variable Lebesgue spaces for the Bergman projector] {Weighted norm inequalities in the variable Lebesgue spaces for the Bergman projector on the unit ball of $\mathbb C^n$}



	\author[D. B\'ekoll\`e]{David B\'ekoll\`e}
	\address{Department of Mathematics, Faculty of Science, University of Yaound\'e I, P.O. Box 812, Yaound\'e, Cameroon }
	\email{{\tt dbekolle@gmail.com}}
	
	\author[E. L. Tchoundja]{Edgar Landry Tchoundja}
\address{Edgar Tchoundja, Leibniz University Hanover\\ Institut f\"ur Analysis
	\\Welfengarden 1 
	\\30167 Hannover, Germany.}
\email{etchoundja@math.uni-hannover.de}
\address{Permanent Address:  Department of Mathematics, Faculty of Science, University of Yaounde I
	\\PO. Box 812
	\\Yaounde-Cameroon. 
}
\email{tchoundjaedgar@yahoo.fr}

\author[A. B. Zotsa-Ngoufack]{Ars\`ene Brice Zotsa Ngoufack}
\address{Aix Marseille Univ, CNRS, I2M, Marseille, France }
\email{\tt arsene-brice.zotsa-ngoufack@univ-amu.fr}

\begin{document}

		\begin{abstract}
			In this work, we  extend the theory of B\'ekoll\`e-Bonami $B_p$ weights. Here we replace the constant $p$ by a non-negative measurable function $p(\cdot),$ which is log-H\"older continuous function with lower bound $1$. We show that the Bergman projector on the unit ball of $\mathbb C^n$ is continuous on the weighted variable Lebesgue spaces $L^{p(\cdot)}(w)$ if and only if $w$ belongs to the generalised B\'ekoll\`e-Bonami class $B_{p(\cdot)}$. To achieve this, we define a maximal function and show that it is bounded on $L^{p(\cdot)}(w)$ if $w\in B_{p(\cdot)}$. We next state and prove a weighted extrapolation theorem that allows us to conclude.
		\end{abstract}
		
		\maketitle

		\section{Introduction}
		The purpose of this work is to generalise the B\'ekoll\`e-Bonami theorem \cite{ref3} for the Bergman projector on the unit ball $\mathbb B$ of $\mathbb C^n$ to the case of weighted variable Lebesgue spaces. The case of the unweighted variable Lebesgue spaces was treated by Chacon and Rafeiro \cite{ref17, ref4}. These authors showed that the Bergman projector is bounded on variable Lebesgue spaces for exponent functions $p(\cdot),$ which are log-H\"older continuous function with lower bound $1$ (cf. Definition \ref{defi1} below).  The ingredients of their proof are: the classical B\'ekoll\`e-Bonami theorem, the boundedness of the Hardy-Littlewood maximal function on variable Lebesgue spaces, and an extrapolation theorem. For basic properties of variable Lebesgue spaces, e.g. the boundedness of Hardy-Littlewood maximal function and an extrapolation theorem, we refer  to \cite{ref16, ref1}.
		
		The $\sigma$-algebra on $\mathbb B$ is the Borel $\sigma$-algebra.	Let $\nu$ be a positive measure on $\mathbb B.$ The variable Lebesgue space on $\mathbb B$, denoted by $L^{p(\cdot)}(\nu)$, is a generalisation of the classical Lebesgue spaces, obtained by replacing the constant exponent $p$ by a measurable exponent function $p(\cdot):\mathbb B\to [0, \infty).$ We shall denote by $\mathcal{P}(\mathbb B)$ this family of all exponent functions $p(\cdot)$ on $\mathbb B.$ For a measurable subset $E$ of $\mathbb B,$ we introduce the following notation:  
		$$p_{-}(E)=ess\inf \limits_{z\in E}p(z)\quad {\rm{ and }}\quad p_{+}(E)=ess\sup \limits_{z\in E}p(z)$$
		and we will use the notation $p_{-}=p_{-}(\mathbb B)$ and $p_{+}=p_{+}(\mathbb B).$  We shall denote  by $\mathcal P_+(\mathbb B),$ the subfamily of $\mathcal{P}(\mathbb B)$ consisting of those $p(\cdot)$ such that $p_+<\infty.$    More precisely, for $p(\cdot)\in \mathcal{P}(\mathbb B),$ we say that $f\in L^{p(\cdot)}(\nu)$ if for some $\lambda>0,\,\rho_{p(\cdot)}(\frac{f}{\lambda})<\infty,$ where 
		$$\rho_{p(\cdot)}(f):=\int_{\mathbb B}|f(z)|^{p(z)}d\nu(z).$$
		For $p(\cdot)\in \mathcal P_+(\mathbb B),$ this definition can be simplified as follows: $f\in L^{p(\cdot)}(\nu)$ if 
		$\rho_{p(\cdot)}(f)<\infty.$
		When $\nu$ is a $\sigma$-finite measure and $p(\cdot)\in \mathcal{P}(\mathbb B)$ is such that $p(\cdot)\geq 1,$ the functional 
		$$\|f\|_{p(\cdot)}=\inf \left\{\lambda>0:\quad \rho_{p(\cdot)}\left (\frac{f}{\lambda}\right )\leq 1\right\}$$
		is a norm on Here, for $z=(z_1,\cdots, z_n)$ and $\zeta=(\zeta_1,\cdots, \zeta_n)\in \mathbb{B}$,  we have set
		$$\langle z,\zeta\rangle=z_1\overline{\zeta_1}+ \cdots +z_n\overline{\zeta_n}\quad\quad \textrm{and} \quad\quad |z|= \langle z,z\rangle^{\frac 12}.$$
		the space $L^{p(\cdot)}(\nu);$ equipped with this norm, $L^{p(\cdot)}(\nu)$ is a Banach space.
		
		We shall denote by $\mu$ the Lebesgue measure on $\mathbb B.$ A non-negative locally integrable function on $\mathbb B$ is called a weight. If $d\nu=wd\mu$ for a weight $w,$  we call $L^{p(\cdot)}(\nu)$ a weighted variable Lebesgue space. In the sequel, $\alpha$ is a positive number and we set $d\mu_{\alpha}(z)=(1-|\zeta|^2)^{\alpha-1}d\mu(z).$ We shall focus on the weighted variable Lebesgue space  $L^{p(\cdot)}(wd\mu_{\alpha}),$ which we shall simply denote $L^{p(\cdot)}(w).$
		
		
		In this paper, we shall take \[d(z,\zeta)=\begin{cases}
		\left \vert |z|-|\zeta|\right \vert+\left \vert 1-\frac{\langle z,\zeta\rangle}{|z||\zeta|}\right \vert&\mbox{if }z,\zeta\in\overline{\mathbb{B}}\setminus\{0\}\\
		|z|+|\zeta|&\mbox{if }z=0\text{ or }\zeta=0
		\end{cases}.\]
		
		This application $d$ is a pseudo-distance on $\mathbb B.$ Explicitly, for all $z,\zeta,\xi\in\mathbb{B},$ we have
		$d(z,\zeta)\leq 2\left (d(z,\xi)+d(\xi,\zeta)\right )$
		and $0\leq d(z, \zeta)<3.$
		In addition, for $z\in\mathbb{B}\text{ and }r>0,$ we denote by \[B(z,r)=\{\zeta\in\mathbb B:\,d(z,\zeta)<r\}\] the open pseudo-ball  centred at $z$ and of radius $r>0$.
		\begin{defn}
			\label{defi1}
			A function $p(\cdot)\in \mathcal P(\mathbb B)$ is  log-H\"older continuous on $\mathbb B$ if there exists a positive constant $c>0$  such that for all $z,\zeta\in\mathbb B$
			\begin{equation}|p(z)-p(\zeta)|\leq\frac{c}{\ln(e+\frac{1}{d(z,\zeta)})}\quad {\rm if} \quad {z\neq \zeta}.\end{equation} 
			We  denote by $\mathcal{P}^{\log}(\mathbb B)$ the space of all log-H\"older continuous functions on $\mathbb B.$ It is easily checked that $\mathcal{P}^{\log}(\mathbb B)\subset \mathcal P_+ (\mathbb B).$ As usual, we set 
			$\mathcal{P}^{\log}_{\pm}(\mathbb B)=\left\{p(\cdot)\in\mathcal{P}^{\log}(\mathbb B):\;p_->1\right\}$.
		\end{defn}
An example of a member of $\mathcal{P}^{\log}_{\pm}(\mathbb B)$ is $p(z)=2+\sin |z|.$
We shall denote by $\mathcal P_-(\mathbb B)$ the subfamily of $\mathcal P(\mathbb B)$ consisting of those $p(\cdot)$ such that $p_->1.$ So $\mathcal{P}^{\log}_{\pm}(\mathbb B)=\mathcal{P}^{\log}(\mathbb B)\cap \mathcal P_-(\mathbb B).$		
		
		\begin{defn}
			We denote by $\mathcal{B}$ the collection of pseudo-balls $B$ of $\mathbb B$ such that $\overline{B}\cap\partial\mathbb B\neq\emptyset$.
			Next we define the maximal function $m_\alpha$ by
			\[m_{\alpha}f(z)=\sup_{B\in\mathcal{B}}\frac{\chi_B(z)}{\mu_{\alpha}(B)}\int_B|f(\zeta)|d\mu_{\alpha}(\zeta).\]
			
			Observe that, from Lemma \ref{lee1}, $\mathcal{B}$ is the set of pseudo-balls that touch the boundary of $\mathbb B.$

		\end{defn}

		In the classical Lebesgue spaces, we have the following definition.
		\begin{defn}
			Let $p>1$ be a constant exponent. The B\'ekoll\`e-Bonami $B_p$ weight class consists of  weights $w$ such that
			\begin{equation*}
			\sup_{B\in\mathcal{B}}\left(\frac{1}{\mu_\alpha(B)}\int_B wd\mu_\alpha\right)\left(\frac{1}{\mu_\alpha(B)}\int_B w^{-\frac{1}{p-1}}d\mu_\alpha\right)^{p-1}<\infty.\label{n-eq-1}
			\end{equation*}
		\end{defn}
		This definition of $B_p$ is equivalent to the following definition:
		$$ \sup_{B\in\mathcal{B}}\frac{1}{\mu_\alpha(B)}\|w^{\frac{1}{p}}\chi_B\|_p\|w^{-\frac{1}{p}}\chi_B\|_{p'}<\infty,$$
		where $p'$ is the conjugate exponent of $p,$ i.e. $\frac 1p+\frac 1{p'}=1.$\\

		In the same spirit, we introduce a variable generalisation  of the $B_p$ weight class. Analogously to the classical case, for $p(\cdot)\in \mathcal P_-(\mathbb B),$ we say that $p'(\cdot)$ is the conjugate exponent function of $p(\cdot)$ if for all $z\in\mathbb B$ we have
		$$\frac{1}{p(z)}+\frac{1}{p'(z)}=1.$$
		Moreover, we set $\mathcal P_\pm (\mathbb B)=\mathcal P_-(\mathbb B)\cap \mathcal P_+(\mathbb B),$  the subfamily of $\mathcal P(\mathbb B)$ consisting of those exponent functions $p(\cdot)$ such that $1<p_{-}\leq p_{+}<\infty.$ We now define the variable B\'ekoll\`e-Bonami classes of weights. 
		\begin{defn}\label{defi2}
			Let $p(\cdot)\in\mathcal{P}_\pm (\mathbb{B}).$  A weight $w$ belongs to the variable B\'ekoll\`e-Bonami class on $\mathbb B$, denoted $B_{p(\cdot)},$ if  
			\begin{equation}
			[w]_{B_{p(\cdot)}}:=\sup_{B\in\mathcal{B}}\frac{1}{\mu_\alpha(B)}\|w^{\frac{1}{p(\cdot)}}\chi_{B}\|_{p(\cdot)}\|w^{-\frac{1}{p(\cdot)}}\chi_{B}\|_{p'(\cdot)}<\infty.\label{neq2}
			\end{equation}
		\end{defn}
		We define the operator $P_{\alpha}$ on $L^1(\mathbb B, d\mu_\alpha)$ by   
		\[P_{\alpha}f(z)=\int_{\mathbb B}\frac{f(\zeta)}{(1-\langle z, \zeta\rangle)^{n+\alpha}}d\mu_{\alpha}(\zeta).\]
		The restriction to $L^2(\mathbb B, d\mu_\alpha)$ of the operator  $P_{\alpha}$ is called the Bergman projector of $\mathbb B$. 
		We also define the positive Bergman operator $P_{\alpha}^+$  by
		\[P_{\alpha}^+ f(z)=\int_{\mathbb B}\frac{f(\zeta)}{\vert 1- \langle z, \zeta\rangle\vert^{n+\alpha}}d\mu_{\alpha}(\zeta).\]
		
		We now recall the classical B\'ekoll\`e-Bonami  theorem: 
		\begin{thm}\cite{ref3}\label{Bek_class}
			Let $w$ be a non-negative measurable function and let $1<p<\infty$ {\rm {($p$ is a constant exponent).}} The following two assertions are equivalent.
\begin{enumerate}
\item[1.]
The Bergman operator $P_\alpha$ is well defined and bounded on $L^p(wd\mu_\alpha);$ 
\item[2.]
$w\in B_p$.
\end{enumerate} 
Moreover, $P_{\alpha}^+$ is well defined and bounded on $L^p(wd\mu_\alpha)$ if $w\in B_p$.
		\end{thm} 
		The purpose of this work is to prove the following generalisation of the previous theorem:
		\begin{thm}\label{th1}
			Let $w$ be a non-negative measurable function and $p(\cdot)\in\mathcal{P}^{\log}_{\pm}(\mathbb B)$. 
The following two assertions are equivalent.
\begin{enumerate}
\item[1.]
The Bergman operator $P_\alpha$ is well defined and bounded on $L^{p(\cdot)}(wd\mu_\alpha);$
\item[2.]
$w\in B_{p(\cdot)}.$
\end{enumerate} 
Moreover, $P_{\alpha}^+$ is well defined and bounded on $L^{p(\cdot)}(wd\mu_\alpha)$ if $w\in B_{p(\cdot)}$. 
		\end{thm} 
The problem under study is trivial if $w(\mathbb B)=0,$ i.e. $w\equiv 0$ a.e. on $\mathbb B.$ We shall assume that $w(\mathbb B)>0.$		
		
In \cite{ref2}, Diening and H\"asto introduced the variable Muckenhoupt weight class $A_{p(\cdot)}$ on $\R^n$, and showed  that for $p(\cdot)\in\mathcal{P}^{\log}_{\pm}(\R^n)$, the Hardy-Littlewood maximal function is bounded on $L^{p(\cdot)}(w)$ only if $w\in A_{p(\cdot)}$. In order to manage the necessary condition, they introduced a new class $A_{p(\cdot)}^+$ which coincides with $A_{p(\cdot)}$ when $p(\cdot)\in\mathcal{P}^{\log}_{\pm}(\R^n),$ but whose condition is easier to check. More precisely, they prove that the class $A_{p(\cdot)}^+$ is contained in the class $A_{p(\cdot)}$ when $p(\cdot)\in\mathcal{P}^{\log}_{\pm}(\R^n),$ but  they leave for future investigation the proof of the reverse inclusion. Later, Cruz-Uribe and these two authors \cite{ref13} gave a new proof of this result using the  Calder\'on-Zygmund decomposition  and they also proved the reverse implication. Very recently, Cruz-Uribe and Cummings \cite{david2022weighted} extended the result of \cite{ref13}  to the spaces of homogeneous type.
		
		In this paper we will use the technique of \cite{ref2, ref13} to manage the proof of the necessary condition in Theorem~\ref{th1}. Precisely, we shall introduce a new class denoted $B_{p(\cdot)}^+$ which coincides with $B_{p(\cdot)}$ when $p(\cdot)\in\mathcal{P}^{\log}_{\pm}(\mathbb B)$. To deal with the sufficient condition, we rely on the result of \cite{david2022weighted} about the boundedness on $L^{p(\cdot)}(w)$ of the Hardy-Littlewood maximal function on the space of homogeneous type $\mathbb B,$ for $p(\cdot)\in\mathcal{P}^{\log}_{\pm}(\mathbb B)$ and $w$ in the corresponding class of weights $A_{p(\cdot)}.$ The proof then follows two steps. First, we use this result to show that, for $p(\cdot)\in\mathcal{P}^{\log}_{\pm}(\mathbb B),$ the maximal function $m_\alpha$ is bounded on $L^{p(\cdot)}(w)$ if $w\in B_{p(\cdot)}.$ Secondly, we lean on the first step to define a new extrapolation theorem which allows us to conclude the proof of the sufficient condition in Theorem~\ref{th1}. 
		
		The rest of our paper is organised as follows. In section~\ref{sec-2}, we shall recall some preliminaries. Next in section~\ref{sec-3}, we review properties of weighted variable Lebesgue spaces, variable B\'ekoll\`e-Bonami and Muckenhoupt classes of weights. In the end of this section, we state the theorem of Cruz-Uribe and Cummings about the boundedness on $L^{p(\cdot)}(w)$ of the Hardy-Littlewood maximal function on $\mathbb{B}.$ In section~\ref{sec-4}, we prove the necessity of the conditions $w^{\frac{1}{p(\cdot)}}\in L^{p(\cdot)}(d\mu_\alpha)$ and $w^{-\frac{1}{p(\cdot)}}\in L^{p'(\cdot)}(d\mu_\alpha)$ in Theorem~\ref{th1}. In section~\ref{sec-5}, we define and study the class $B^+_{p(\cdot)}$ and we show the identity $B^+_{p(\cdot)}=B_{p(\cdot)}$ and we also prove the reverse inclusion $A_{p(\cdot)}$ is contained in $A_{p(\cdot)}^+$. In section~\ref{sec-6}, we prove the necessary condition in Theorem~\ref{th1}. In section~\ref{sec-7}, we show that the maximal function $m_\alpha$ is bounded on $L^{p(\cdot)}(w)$ if $w\in B_{p(\cdot)}.$ Finally, in section~\ref{sec-8}, we prove a weighted extrapolation theorem from which we deduce a proof of the sufficient condition in Theorem~\ref{th1}.
		
		Let $a$ and $b$ be two positive numbers. Throughout the paper, we write $a\lesssim b$ if there exists $C>0$ such that $a\leq Cb$. We write $a\simeq b$ if $a\lesssim b$ and $b\lesssim a$.

		\section{Preliminaries} \label{sec-2}
		In this section we present some background material regarding the unit ball of $\mathbb{C}^n$ as a space of homogeneous type and the variable exponent Lebesgue spaces. 
		\subsection{The unit ball  is a space of homogeneous type}
		
		In this subsection we recall some lemmas from \cite{ref3} where it was first considered $(\mathbb B, d, \mu_\alpha)$.
		
		\begin{lem}\label{lee1}
			Let $z\in\mathbb B$ and $R>0$. The pseudo-ball $B(z,r)$ meets the boundary of $\mathbb B$ (in other words, $B\in \mathcal B$) if and only if $R> 1-|z|$.
		\end{lem}
		\begin{lem}\label{lee2}
			Let $z\in\mathbb B$ and $0<R\leq3$, we have 
			\[\mu_{\alpha}(B(z,R))\simeq R^{n+1}\left (\max(R,1-|z|)\right )^{\alpha-1}.\]
		\end{lem}
		\begin{rem}\label{rqq1}
			From Lemma~\ref{lee1} and Lemma~\ref{lee2} if $B(z,R)\in\mathcal{B}$ we have \[\mu_{\alpha}(B(z,R))\simeq R^{n+\alpha}.\]
		\end{rem}
		From Lemma~\ref{lee2},  that $(\mathbb B,d,\mu_{\alpha})$ is a space of homogeneous type. Next we have the following lemma.
		\begin{lem}\cite{david2022weighted}\label{lem12}
			There exist two positive constants $C$ and $\gamma$ such that, for all $\zeta\in B(z,R)$ we have
			\[\mu_{\alpha}(B(\zeta,r))\geq C\left(\frac{r}{R}\right)^{\gamma}\mu_{\alpha}(B(z,R))\] for every $0< r\leq R<\infty.$
		\end{lem}
		
		\subsection{Variable exponent Lebesgue spaces}
		We denote by $\mathcal{M}$ the space of complex-valued measurable functions defined on $\mathbb B$. Let $\nu$ be a positive measure on $\mathbb B.$ The family $\mathcal{P}(\mathbb B)$ of variable exponents was defined in the introduction. 
		In the rest of the paper, we take $p(\cdot)\in \mathcal{P}(\mathbb B).$
		 The next definitions, properties and propositions are stated in  \cite{ref16,ref1,ref4}. We first recall some properties of the modular functional $\rho_{p(\cdot)}:\mathcal M\rightarrow [0, \infty],$ defined in the introduction as 
		$$\rho_{p(\cdot)} (f)=\int_{\mathbb B} |f(z)|^{p(z)}d\nu (z).$$
		
		\begin{prop}
			\label{prop-rho}
			Let $p(\cdot)\in\mathcal{P}(\mathbb B)$ be such that $p(\cdot)\geq 1.$
			\begin{enumerate}
				\item For all $f\in\mathcal{M},\quad\rho_{p(\cdot)}(f)\geq 0 ,\text{ and }\rho_{p(\cdot)}(f)=\rho_{p(\cdot)}(|f|)$.
				\item  For all $f\in\mathcal{M}\text{ if } \rho_{p(\cdot)}(f)<\infty$ then $|f(z)|<\infty\text{ a.e.}$ on $\mathbb B.$
				\item \label{pro3}$\rho_{p(\cdot)}$ is convex. In particular, for $0<\alpha\leq1$ and $f\in\mathcal{M},\,\rho_{p(\cdot)}(\alpha f)\leq\alpha\rho_{p(\cdot)}(f)$ and for $\alpha\geq 1,\quad \alpha\rho_{p(\cdot)}(f)\leq\rho_{p(.)}(\alpha f)$.	 
				\item $\rho_{p(\cdot)}(f)=0 \text{ if and only if } f(z)=0$ a.e. on $\mathbb B.$
				\item If for almost all $z\in\mathbb B,\,|f(z)|\leq|g(z)|$, then $\rho_{p(\cdot)}(f)\leq\rho_{p(\cdot)}(g)$.
				\item \label{pro7} If there exists $\beta>0$ such that $\,\rho_{p(\cdot)}(\frac{f}{\beta})<\infty$, then the function  $\lambda\longmapsto\rho_{p(\cdot)}(\frac{f}{\lambda})$ is continuous and non-increasing on $[\beta,\infty[$. In addition \[\lim_{\lambda\longrightarrow \infty}\rho_{p(.)}\left (\frac{f}{\lambda}\right )=0\].  
			\end{enumerate}
		\end{prop}

		For $p(\cdot)\in\mathcal{P}(\mathbb B),$ the variable Lebesgue space $L^{p(\cdot)}(d\nu)$ was defined in the introduction.
		
		\begin{prop}\cite[Theorem 2.7.2]{ref16}\label{dense}
			Let $p(\cdot)\in\mathcal{P}_+(\mathbb B)$ be such that $p(\cdot)\geq 1.$  Then the subspace of continuous functions of compact support in $\mathbb B$ is dense in the space $L^{p(\cdot)}(d\nu).$ 
		\end{prop}
		
		We next recall the H\"older inequality in the variable exponent context.
		
		\begin{prop}
			\cite[Theorem 2.26, Corollary 2.28]{ref16}\label{hold} 
			\begin{enumerate}
				\item[1.]
				Let $p(\cdot)\in\mathcal{P}(\mathbb B)$ be such that $p(\cdot)\geq 1.$  Then for all $f,g\in \mathcal M,$ we have
				\[\int_{\mathbb B}|fg|d\nu\leq2\|f\|_{p(\cdot)}\|g\|_{p'(\cdot)}.\]
				\item[2.]
				Let $r(\cdot), \hskip 1truemm q(\cdot)\in \mathcal P(\mathbb B)$ such that $r(\cdot), q(\cdot)\geq 1$ and $\frac 1{q(x)}+\frac 1{r(x)}\leq 1$ for all $x\in \mathbb B.$ Define $p(\cdot)\in \mathcal P(\mathbb B)$ such that $p(\cdot)\geq 1,$ by
				$$\frac 1{p(x)}=\frac 1{q(x)}+\frac 1{r(x)}.$$
				Then there exists a constant $K$ such that for all $f\in L^{q(\cdot)}$ and $g\in L^{r(\cdot)},$ $fg\in  L^{p(\cdot)}$ and 
				$$\Vert fg\Vert_{p(\cdot)}\leq K\|f\|_{q(\cdot)}\|g\|_{r(\cdot)}.$$ 
			\end{enumerate} 
		\end{prop}
We record the following useful remark.
\begin{rem}\label{p-p'}
			 The property   
			$p(\cdot)\in \mathcal{P}^{\log}_{\pm}(\mathbb B)$ is also true for $p'(\cdot).$ 
		\end{rem}

		In what follows, we will use condition log-H\"older as in given in the following Lemma.
		
		\begin{lem}\label{log}
			Let $p(\cdot)\in\mathcal{P}^{\log}(\mathbb B).$  Let $B=B(x, R)$ be a pseudo-ball of $\mathbb B$ such that $R<\frac 14.$ Then 
			\[p_{+}(B)-p_{-}(B)\leq\frac{c}{\ln(\frac 1{4R})}.\]
		\end{lem}

		\begin{lem}\label{lee3}
			Let $p(\cdot)\in\mathcal{P}^{\log}(\mathbb B)$. There exist two positive constants $C_1=C_1(\alpha,n,p(\cdot))$ and $C_2=C_2(\alpha,n,p(\cdot))$ such that for every pseudo-ball $B$ of $\mathbb B,$ we have
			\[\mu_{\alpha}(B)^{p_{-}(B)-p_{+}(B)}\leq C_1\text{ and }\mu_{\alpha}(B)^{p_{+}(B)-p_{-}(B)}\leq C_2.\]
		\end{lem}
		\begin{proof}
			Since $p_{+}(B)-p_{-}(B)\geq0,\,\mu_{\alpha}(\mathbb B)<\infty,\,B\subset\mathbb B$ and $p(\cdot)$ bounded, we have the second inequality. 
   We are going to prove the first inequality. Suppose that $B=B(z,R)$.
			\begin{enumerate}
				\item[1.] If $R\geq\frac{1}{16}$, from Lemma~\ref{lee2}, there exists $C>0$ such that 
				$\mu_{\alpha}(B)\geq C R^{n+\alpha}$
				and as  $p_{-}(B)-p_{+}(B)\leq 0$, we obtain:
				\begin{equation*}
				\mu_{\alpha}(B)^{p_{-}(B)-p_{+}(B)}\leq (CR^{n+\alpha})^{(p_{-}(B)-p_{+}(B))}\leq (C(16)^{n+1})^{(p_{+}(B)-p_{-}(B))}
				\end{equation*}
				\item[2.] If $R<\frac{1}{16}$, 
				from Lemma \ref{log}, 
				\[p_{+}(B)-p_{-}(B)\leq\frac{c}{\ln(\frac 1{4R})}.\]
				Hence from Remark \ref{rqq1},
				\[p_{+}(B)-p_{-}(B)\lesssim\frac{1}{\ln(\frac 1{\mu_{\alpha}(B)})}.\]
				Thus 
				\begin{equation*}
					\mu_{\alpha}(B)^{p_{-}(B)-p_{+}(B)}\lesssim1.
				\end{equation*}
			\end{enumerate}
		\end{proof}
		
		From the previous lemma, we easily deduce the following corollary.
		\begin{cor}\label{rqq2}
			Let $p(\cdot)\in\mathcal{P}^{\log}(\mathbb B).$  There exists a  constant $C=C(\alpha,n,p(\cdot))>1$  such that for every pseudo-ball $B$ of $\mathbb B$ and every $z\in B$ we have
			\[\frac 1C\leq \mu_{\alpha}(B)^{p_{-}(B)-p(z)}\leq C.\]
		\end{cor}
		
		\section{Weighted variable  Lebesgue spaces,  variable B\'ekoll\`e-Bonami  and Muckenhoupt classes of weights}\label{sec-3}
\subsection{Weighted variable  Lebesgue spaces}
		Let $w$ be a weight and let $p(\cdot)\in\mathcal{P}(\mathbb B)$ be such that $p(\cdot)\geq 1.$ The corresponding weighted variable exponent Lebesgue space $L^{p(\cdot)}(wd\mu_{\alpha})$ consists of those $f\in \mathcal M$ which satisfy the estimate
		\[\rho_{p(\cdot),w}(f)=\int_{\mathbb B}|f(z)|^{p(z)}w(z)d\mu_{\alpha}(z)<\infty.\]
		We also denote it $L^{p(\cdot)}(\mathbb B,w)$ or simply $L^{p(\cdot)}(w),$ and we denote its norm by $\Vert \cdot\Vert_{p(\cdot),w}.$ It is easy to check that 
		\[\|f\|_{p(\cdot),w}=\|fw^{\frac{1}{p(\cdot)}}\|_{p(\cdot)}.\]
		In the sequel, we shall adopt the following notation:
		\begin{equation}\label{eqomegaprime}
		      w':=w^{1-p'(\cdot)}.
		\end{equation}
 
		We recall the notion of subordinate norm on $L^{p(\cdot)}(\mathbb B,w)$ defined by
		\[\Vert f\Vert'_{p(\cdot), w}:=\sup_{\|g\|_{p'(\cdot),w'}=1}\left|\int_{\mathbb B}f(\zeta)\overline{g(\zeta)}d\mu_{\alpha}(\zeta)\right|.\]
		We next recall the following proposition. 
		\begin{prop}\cite[Corollary $2.7.5$]{ref16,ref1} \label{egalm}
			Let $p(\cdot)\in\mathcal{P}(\mathbb B)$ such that $p(\cdot)\geq 1$ and let $w$ be a weight. Then
			\[\|f\|_{p(\cdot),w}\leq\Vert f\Vert'_{p(\cdot), w} \leq2\|f\|_{p(\cdot),w}.\]
		\end{prop}
		The following lemma will be very useful.
		\begin{lem}\cite{ref16, ref1}\label{le8}
			Let $w$ be a non-negative measurable function and let $p(\cdot)\in\mathcal{P_+}(\mathbb B)$ be such that $p_->0$. Then for every $f\in \mathcal M$ whose support is $E,$ the following double inequality holds 
			\begin{equation*}
			\min\left(\rho_{p(\cdot),w}(f)^{\frac{1}{p_-(E)}},\rho_{p(\cdot),w}(f)^{\frac{1}{p_+(E)}}\right)
			\leq\|f\|_{p(\cdot),w}\leq\max\left(\rho_{p(\cdot),w}(f)^{\frac{1}{p_-(E)}},\rho_{p(\cdot),w}(f)^{\frac{1}{p_+(E)}}\right).
			\end{equation*} 
			It is equivalent to
			\begin{equation*}
			\min\left(\|f\|_{p(\cdot),w}^{p_-(E)},\|f\|_{p(\cdot),w}^{p_+(E)}\right)
			\leq\rho_{p(\cdot),w}(f)\leq\max\left(\|f\|_{p(\cdot),w}^{p_-(E)},\|f\|_{p(\cdot),w}^{p_+(E)}\right).
			\end{equation*}
		\end{lem}

\subsection{Variable B\'ekoll\`e-Bonami classes of weights}		
		Concerning the variable B\'ekoll\`e-Bonami weight class  $B_{p(\cdot)},$ we record the following elementary remark. 

\begin{rem}\label{elem}
Let $p(\cdot)\in \mathcal P_\pm(\mathbb B).$ If $w\in B_{p(\cdot)},$ the following two assertions are valid.
\begin{enumerate}
\item[1)]
$\Vert w^{\frac {1}{p(\cdot)}}\Vert_{p(\cdot)}<\infty$ and $\Vert w^{-\frac {1}{p(\cdot)}}\Vert_{p'(\cdot)}<\infty.$
\item[2)]
The functions $w$ and $w'$ are integrable on $\mathbb B.$
\end{enumerate}
\end{rem}

\begin{proof}[Proof of the remark]
Otherwise, if $\Vert w^{\frac {1}{p(\cdot)}}\Vert_{p(\cdot)}=\infty,$ then necessarily $\Vert w^{-\frac {1}{p(\cdot)}}\Vert_{p'(\cdot)}=0$ and this would imply that $w\equiv \infty$ a.e. Alternatively, if $\Vert w^{-\frac {1}{p(\cdot)}}\Vert_{p(\cdot)}=\infty,$ then necessarily $\Vert w^{\frac {1}{p(\cdot)}}\Vert_{p'(\cdot)}=0$ and this would imply that $w\equiv 0$ a.e. Furthermore, by Lemma~\ref{le8}, for $p(\cdot)\in \mathcal P_\pm(\mathbb B),$  the assertions 1) and 2) are equivalent.
\end{proof}

We also have the following proposition.
		
		\begin{prop}\label{w-w'}
			Let $p(\cdot)\in \mathcal {P}_\pm(\mathbb B).$ For a weight $w,$  the following two assertions are equivalent.
			\begin{enumerate}
				\item
				$w\in B_{p(\cdot)};$
				\item
				$w'\in B_{p'(\cdot)}.$
			\end{enumerate}
			Moreover, $[w]_{B_{p(\cdot)}}=[w']_{B_{p'(\cdot)}}.$
		\end{prop}
The following simple lemma will be useful.
\begin{lem}
Let $p(\cdot)\in \mathcal {P}_\pm(\mathbb B).$ For a weight $w,$  the following two assertions are equivalent.
			\begin{enumerate}
				\item
				$w\in B_{p(\cdot)};$
				\item
$\sup \limits_{B\in \mathcal B}\frac 1{\mu_\alpha (B)}\left \Vert \chi_B\right \Vert_{p(\cdot), w}\left \Vert \chi_B\right \Vert_{p'(\cdot), w'}<\infty.$
       \end{enumerate} 
\end{lem}

\subsection{Variable Muckenhoupt classes of weights}
		
		\begin{defn}
			The Hardy-Littlewood maximal function $M_\alpha$ on the space of homogeneous type $(\mathbb B, d, \mu_\alpha)$ is defined by
			\[M_{\alpha}f(z)=\sup_{B}\frac{\chi_B(z)}{\mu_{\alpha}(B)}\int_B|f(\zeta)|d\mu_{\alpha}(\zeta)\] where the supremum is taken over all pseudo-balls of $\mathbb B.$
		\end{defn}
		
		When $p$ is a constant greater than $1$,  the Muckenhoupt class $A_p$ consists of weights $w$ which satisfy the estimate
		\[\sup_{B}\left(\frac{1}{\mu_\alpha(B)}\int_B wd\mu_\alpha\right)\left(\frac{1}{\mu_\alpha(B)}\int_B w^{-\frac{1}{p-1}}d\mu_\alpha\right)^{p-1}<\infty,\]
		where the $\sup$ is taken over all pseudo-balls $B$ of $\mathbb B.$ This definition is equivalent with the following definition:
		\[\sup_{B}\frac{1}{\mu_\alpha(B)}\|w^{\frac{1}{p}}\chi_B\|_p\|w^{-\frac{1}{p}}\chi_B\|_{p'}<\infty,\] 
		where the $\sup$ is taken again over all pseudo-balls of $\mathbb B.$
		
		We next have the following variable generalisation  of the variable Muckenhoupt weight classes. This generalisation was given first by Diening and H\"asto \cite{ref2}. 
		\begin{defn}
			Let $p(\cdot)\in\mathcal{P}_\pm(\mathbb{B}).$  A weight $w$ belongs to the variable Muckenhoupt class $A_{p(\cdot)}$ on $\mathbb B$ if  
			\begin{equation}
			[w]_{A_{p(\cdot)}}:=\sup_{B}\frac{1}{\mu_\alpha(B)}\|w^{\frac{1}{p(\cdot)}}\chi_{B}\|_{p(\cdot)}\|w^{-\frac{1}{p(\cdot)}}\chi_{B}\|_{p'(\cdot)}<\infty\label{eq1}
			\end{equation}
			where the $\sup$ is taken over all pseudo-balls of $\mathbb B.$
			
		\end{defn}
		Let $p(\cdot)\in\mathcal{P}_\pm(\mathbb{B}).$ The following proposition is similar to Proposition \ref{w-w'}.
		\begin{prop}\label{w-w'bis}
			The following two assertions are equivalent.
			\begin{enumerate}
				\item
				$w\in A_{p(\cdot)};$
				\item
				$w'\in A_{p'(\cdot)}.$
			\end{enumerate}
		\end{prop}
		We record the following properties of $A_{p(\cdot)}$ and $B_{p(\cdot)}.$
		\begin{prop}\label{rq1}
			\begin{enumerate}
				\item[1.] The inclusion $A_{p(\cdot)}\subset B_{p(\cdot)}$ holds with $[w]_{B_{p(\cdot)}}\leq [w]_{A_{p(\cdot)}}.$
				\item[2.] $[w]_{B_{p(\cdot)}}\geq \frac 12$ and $[w]_{A_{p(\cdot)}}\geq \frac 12$.
			\end{enumerate}
		\end{prop}
		
		\begin{proof}
			\begin{enumerate}
				\item[1.]
				This follows directly from the definitions of $A_{p(\cdot)}$ and $B_{p(\cdot)}.$
				\item[2.]		
				We first give the proof for $B_{p(\cdot)}.$ Let $B\in \mathcal B.$		From the H\"older inequality and the definition of $B_{p(\cdot)},$ we have
				\begin{align}
				1=\frac{1}{\mu_\alpha(B)}\int_Bw^{\frac{1}{p(\cdot)}}w^{-\frac{1}{p(\cdot)}}d\mu_\alpha\leq \frac{2}{\mu_\alpha(B)}\|w^{\frac{1}{p(\cdot)}}\chi_{B}\|_{p(\cdot)}\|w^{-\frac{1}{p(\cdot)}}\chi_{B}\|_{p'(\cdot)}\leq 2[w]_{B_{p(\cdot)}}.
				\end{align}
				The proof for $A_{p(\cdot)}$ then follows from assertion 1.
			\end{enumerate}
		\end{proof}
		In \cite{david2022weighted}, Cruz-Uribe and Cummings proved the following fundamental result for the maximal Hardy-Littlewood function. This variable theorem generalises a well-known theorem of Muckenhoupt \cite{M72} in the Euclidean space $\mathbb R^n$. For spaces of homogeneous type, the analogous theorem for constant exponents was proved later by A. P. Calder\'on \cite{C76}. 
		\begin{thm}\label{max}
			Let $p(\cdot)\in\mathcal{P}^{\log}_{\pm}(\mathbb B)$. The following two assertions are equivalent.
			\begin{enumerate}
				\item[1.] 
				There exists a positive constant $C$ such that for all $f\in L^{p(\cdot)}(w),$ we have
				\[\|M_\alpha f\|_{p(\cdot),w}\leq C\|f\|_{p(\cdot),w}.\]
				\item[2.]
				$w\in A_{p(\cdot)}.$
			\end{enumerate} 
		\end{thm}

In fact, these authors \cite{david2022weighted} proved their result in the general setting of spaces of homogeneous type. There, in addition to the condition $p(\cdot)\in \mathcal P^{\log} (\mathbb B)$, they have a condition at infinity which has the following expression on the unit ball $\mathbb B$: there are two constants $c$ and $p_\infty$ such that 
\begin{equation*}|p(z)-p_\infty|\leq\frac{c}{\ln(e+|z|)}\end{equation*}
for every $z\in \mathbb B.$ It is easy to show that this extra condition is satisfied on $\mathbb B.$

Since $m_\alpha f\leq M_\alpha f,$ we deduce the following corollary. 
		\begin{cor}\label{rqinv}
			 Let $p(\cdot)\in\mathcal{P}^{\log}_{\pm}(\mathbb B)$ and $w\in A_{p(\cdot)}.$ For the same constant $C$ as in Theorem ~\ref{max}, we have  
			\[\|m_\alpha f\|_{p(\cdot),w}\leq C\|f\|_{p(\cdot),w}\]
			for all $f\in L^{p(\cdot)}(w).$
		\end{cor}

		\section{Necessity of the conditions $w^{\frac{1}{p(\cdot)}}\in L^{p(\cdot)}(d\mu_\alpha)$ and $w^{-\frac{1}{p(\cdot)}}\in L^{p'(\cdot)}(d\mu_\alpha)$ in Theorem~\ref{th1}}\label{sec-4}
		\begin{prop}\label{l1}
			Let $w$ be a weight and  let $p(\cdot)\in\mathcal{P}_- (\mathbb B).$    If the operator $P_{\alpha}$ is well defined on $L^{p(\cdot)}(w),$ then $w^{-\frac{1}{p(\cdot)}}\in L^{p'(\cdot)}(d\mu_\alpha)$.
		\end{prop}
		\begin{proof}
			If the operator $P_{\alpha}$ is well defined on $L^{p(\cdot)}(w),$ then for every $f\in L^{p(\cdot)}(w),$ we have 
			\[P^+_\alpha |f|(z)=\int_{\mathbb{B}}\frac{|f(\zeta)|}{|1-z\cdot\overline{\zeta}|^{n+\alpha}}d\mu_{\alpha}(\zeta)<\infty\quad \text{for all } \hskip 2truemm   \hskip 2truemm z\in\mathbb B.\]
			
			By a contradiction argument, suppose that $w^{-\frac{1}{p(\cdot)}}$ does not belong to $L^{p'(\cdot)}(d\mu_\alpha).$ Then by Proposition 3.1, there exists a non-negative   $g\in L^{p(\cdot)}(d\mu_\alpha)$  such that  \[\int_{\mathbb B}g(\zeta)w(\zeta)^{-\frac{1}{p(\zeta)}}d\mu_{\alpha}(\zeta)=\infty.\]	
			Let $f=gw^{-\frac{1}{p(\cdot)}}.$ We have $f\in L^{p(\cdot)}(w)$ but $f$ does not belong to $L^1(\mathbb{B});$ and as $|1-\overline{\zeta}z|\leq2$ then \[P^+_\alpha |f|(z)=\int_{\mathbb{B}}\frac{|f(\zeta)|}{|1-\overline{\zeta}\cdot z|^{n+\alpha}}d\mu_{\alpha}(\zeta)\geq\frac{1}{2^{n+\alpha}}\int_{\mathbb{B}}|f(\zeta)|d\mu_{\alpha}(\zeta)=\infty.\]
			This contradict the fact that $P^+_\alpha |f|(z)<\infty$ and consequently we have the result.
		\end{proof}
		\begin{prop}\label{l2}
			Let $w$ be a weight and $p(\cdot)\in\mathcal{P}_-(\mathbb B).$  If the operator $P_{\alpha}$ is bounded on $L^{p(\cdot)}(wd\mu)$, then $w^{\frac{1}{p(\cdot)}}\in L^{p(\cdot)}(d\mu_\alpha)$.
		\end{prop}
		\begin{proof}
			Let $0<r<1$ and define the function $f(z)=(1-|z|^2)^{1-\alpha}\chi_{B(0,r)}(z)$ on $\mathbb B.$ 
   We have 
			\begin{align*}
			P_{\alpha}f(z)  &=\int_{\mathbb{B}}\frac{f(\zeta)}{\left(1-\langle z, \zeta\rangle\right)^{n+\alpha}}d\mu_{\alpha}(\zeta)\\
			&=\int_{B(0,r)}\frac{1}{\left(1-\langle z, \zeta\rangle\right)^{n+\alpha}}d\mu(\zeta)\\
			&=\overline{\int_{B(0,r)}\frac{1}{\left(1-\langle \zeta, z \rangle\right)^{n+\alpha}}d\mu(\zeta)}.
			\end{align*}
			Since the function $\zeta\longmapsto\, \frac{1}{\left(1-\langle \zeta, z \rangle\right)^{n+\alpha}}$ is analytic on $\mathbb{B}$ and $B(0, r)$ is the Euclidean ball centred at $0$ and of radius $r,$ it follows from the mean value property that
			$P_{\alpha}f(z)\equiv C_{r,n}$
			and so
\begin{equation}\label{mean}
|C_{r,n}|\|w^{\frac{1}{p(\cdot)}}\|_{p(\cdot)}=\|P_{\alpha}(f)\|_{p(\cdot),w}.
\end{equation}			
			In addition, 
			\begin{align*}
			\rho_{p(\cdot), w}(f)=\int_{B(0,r)}w(z)(1-|z|^2)^{(\alpha-1)(1-p(z))}d\mu(z).
			\end{align*}
			On the one hand, if $\alpha\leq1$, we have  $(1-|z|^2)^{(\alpha-1)(1-p(z))}\leq1$ because $(\alpha-1)(1-p(z))>0$ and $1-|z|^2\leq1$.
			Consequently
			\[\rho_{p(\cdot), w}(f)\leq\int_{B(0,r)}w(z)d\mu(z)<\infty\] because $w$ is locally integrable.
			\\On the other hand, if $\alpha>1$ we have $(1-|z|^2)^{(\alpha-1)(1-p(z))}\leq(1-|z|^2)^{(\alpha-1)(1-p_{+})}.$
			So
			\[\rho_{p(\cdot), w}(f)\leq\sup_{z\in B(0,r)}(1-|z|^2)^{(\alpha-1)(1-p_{+})}\int_{B(0,r)}w(z)d\mu(z)<\infty\]because $w$ is locally integrable and $\sup_{z\in B(0,r)}(1-|z|^2)^{(\alpha-1)(1-p_{+})}=(1-r^2)^{(\alpha-1)(1-p_{+})}.$ 
			\\Thus, since $\rho_{p(\cdot),w}(f)<\infty$ in both cases, by Lemma \ref{le8} we obtain $\|f\|_{p(\cdot),w}<\infty$ and as $P_{\alpha}$ is bounded on $L^{p(\cdot)}(w),$ we deduce from \eqref{mean} that there exists a positive constant $c_{r,\alpha,n}$ such that 
			\[\|w^{\frac{1}{p(\cdot)}}\|_{p(\cdot)}\leq c_{r,\alpha,n}\|f\|_{p(\cdot),w}<\infty.\] Hence we have the result.
		\end{proof}
		In what follows, we need to calculate the $p(\cdot)-$norm of some functions but is not easy like a constant case. To deal with that, in the following section we introduce the definitions of new spaces whose coincide with the class $B_{p(\cdot)}$.    
		\section{The weight classes $B^+_{p(\cdot)}$ and $B^{++}_{p(\cdot)}$}\label{sec-5}
		\begin{defn}\label{defi3}
			Let $p(\cdot)\in\mathcal{P}_{\pm}(\mathbb{B})$ and let $w$ be a weight. We say that $w$ is in the $B^{+}_{p(\cdot)}$ class if 
			\begin{equation}[w]_{B^{+}_{p(\cdot)}}:=\sup \limits_{B\in \mathcal B}\frac 1{\mu_{\alpha}(B)^{p_B}}\|w\chi_{B}\|_1\|w^{-1}\chi_B\|_{\frac {p'(\cdot)}{p(\cdot)}}<\infty.\label{eq2}\end{equation}  
			where 
   $$p_B=\left (\frac{1}{\mu_{\alpha}(B)}\int_B \frac 1{p(x)}d\mu_{\alpha}(x)\right )^{-1}.$$
		\end{defn}
		
		This class coincides with the $B_p$ class when $p(\cdot)=p \hskip 2truemm (p\text{ constant})$.	We also adopt the following notation:
\begin{equation}\label{eq-meanexponent}
  \langle p\rangle_B=\frac{1}{\mu_{\alpha}(B)}\int_B p(x)d\mu_{\alpha}(x).  
\end{equation}

		\begin{rem}\label{rq2}
			Let $p(\cdot)\in\mathcal{P}^{\log}(\mathbb B)$ and let $B$ be a pseudo-ball of $\mathbb B$. As $p_{-}(B)\leq p_{B}, \langle p\rangle_B\leq p_{+}(B)$, it follows from Lemma~\ref{lee3} and Corollary~\ref{rqq2} that  \[\mu_{\alpha}(B)^{p_{-}(B)}\simeq \mu_{\alpha}(B)^{p_{B}}\simeq \mu_{\alpha}(B)^{p_{B}}\simeq \mu_{\alpha}(B)^{\langle p\rangle_{B}} \simeq\mu_{\alpha}(B)^{p_{+}(B)}.\]
		\end{rem}
		\begin{lem}\label{le4}\cite[Theorem 4.5.7]{ref1}\\
			Let $p(\cdot)\in\mathcal{P}^{log}(\mathbb B)$ be such that $p_->0.$ Let $B$ be a pseudo-ball of $\mathbb B$. Then \[\|\chi_{B}\|_{p(\cdot)}\simeq\mu_{\alpha}(B)^{\frac{1}{p_B}}.\]
		\end{lem}
		\begin{proof}
			From Lemma~\ref{le8} we have
			\begin{equation*}
			\min\left(\mu_\alpha(B)^{\frac{1}{p_-(B)}},\mu_\alpha(B)^{\frac{1}{p_+(B)}}\right)
			\leq\|\chi_B\|_{p(\cdot)}\leq\max\left(\mu_\alpha(B)^{\frac{1}{p_-(B)}},\mu_\alpha(B)^{\frac{1}{p_+(B)}}\right).
			\end{equation*}
			Next, from the Remark~\ref{rq2}, we have $\mu_\alpha(B)^{\frac{1}{p_-(B)}}\simeq\mu_\alpha(B)^{\frac{1}{p_+(B)}}\simeq\mu_\alpha(B)^{\frac{1}{p_B}}$ because $p(\cdot)$ is bounded away from zero. The conclusion follows. 
		\end{proof}
		\begin{lem}\label{le5}
			Let $p(\cdot)\in \mathcal P^{log}_\pm (\mathbb B)$  and let $q$ be a constant exponent greater than $p_++1.$ There exists a positive constant $C$ depending only of the log-H\"older constant of $p(\cdot)$ such that
			\[[w]_{B_{q}}\leq C [w]_{B_{p(\cdot)}^{+}}.\]

		\end{lem}
		\begin{proof}
			As $p(\cdot)<p_+ +1<q,$ we have $\frac {q'}q<\frac {p'(\cdot)}{p(\cdot)}.$ Hence from the H\"older inequality (assertion 2 of Proposition \ref{hold}), we obtain
			\[\|w^{-1}\chi_{B}\|_{\frac {q'}{q}}\leq K\|\chi_{B}\|_{\beta(\cdot)}\|w^{-1}\chi_{B}\|_{\frac{p'(\cdot)}{p(\cdot)}},\]
			where 
			\[\frac{1}{\beta(\cdot)}=\frac{q}{q'}-\frac{p(\cdot)}{p'(\cdot)}=q-p(\cdot)>1.\]
			It is easy to check that $\beta(\cdot)$ is a member of $\mathcal P^{log}(\mathbb B)$ such that $\beta_->0.$ Consequently, from Lemma \ref{le4} and Remark \ref{rq2}, we have
			\[\|\chi_{B}\|_{\beta(\cdot)}\simeq \mu_{\alpha}(B)^{\frac{1}{\beta_B}}\simeq\mu_{\alpha}(B)^{q-\langle p \rangle_B}\simeq\mu_{\alpha}(B)^{q- p_B}.\]
			Thus there exists a positive constant $C$ such that
			$$
			\begin{array}{clcr}
			\frac{1}{\mu_{\alpha}(B)^q}\|w\chi_{B}\|_1\|w^{-1}\chi_B\|_{\frac {q'}{q}}&\leq C\frac{1}{\mu_{\alpha}(B)^{q}}\|w\chi_{B}\|_1\|w^{-1}\chi_B\|_{\frac {p'(\cdot)}{p(\cdot)}}\mu_\alpha (B)^{q-p_B}\\
			&=C\frac{1}{\mu_{\alpha}(B)^{p_B}}\|w\chi_{B}\|_1\|w^{-1}\chi_B\|_{\frac {p'(\cdot)}{p(\cdot)}}\leq C[w]_{B^+_{p(\cdot)}}.
			\end{array}
$$
for all pseudo-balls $B\in \mathcal B.$ The conclusion follows.
			
		\end{proof}
		
		We recall the following definition.
		\begin{defn}
			The weight class $B_\infty$ is defined by $B_\infty=\bigcup \limits_{q\in (1,\infty)} B_q.$
		\end{defn}
		
		\begin{rem}\label{embed}
			It follows from Lemma~\ref{le5} that if $p(\cdot)\in \mathcal P^{\log}_\pm (\mathbb B),$  we have
			\[B_{p(\cdot)}^{+}\subset B_{\infty}.\]
		\end{rem}
		In the rest of this article, to simplify the notation, we denote $w(B)=\|w\chi_B\|_1$.
		
		We next define another class of weights $\Lambda$, which contains the class $B_\infty.$ For a reference, cf. e.g.  \cite{APR2019}.
		\begin{defn}\label{def-lambda}
			We call $\Lambda$ the class consisting of those integrable weights $w$ satisfying the following property. There exist two positive constants $C$ and $\delta$  such that the following implication holds.
\begin{equation}\label{eq-lambda}
\frac {\mu_\alpha (E)}{\mu_\alpha (B)}\leq C \left (\frac {w(E)}{w(B)}\right )^\delta
\end{equation}

			whenever $B\in \mathcal{B}$ and $E$ is a measurable subset of $B.$
		\end{defn}
		\begin{rem}\label{doubling}
\begin{enumerate}
\item[1.]
			For $w\in \Lambda,$ the weighted measure $wd\mu_\alpha$ is doubling in the following sense. There exists a positive constant $C$ such that for every pseudo-ball $B$ of $\mathbb B$ whose pseudo-ball $\widetilde B$ of same centre and of double radius  is a member of $\mathcal B$, we have
			$$w\left (\widetilde B \right )\leq Cw(B).$$
			This result easily follows from the definition of $\Lambda.$ 
\item[2.]
We recall that $w(\mathbb B)>0.$ For $w\in \Lambda,$ this implies that $w(B)>0$ for every pseudo-ball $B$ of $\mathbb B.$  Indeed, take $\mathbb B$ for $B$ and $B$ for $E$ in \eqref{eq-lambda}.
\end{enumerate}
		\end{rem}
		
		
		\begin{lem}\label{le9}
			Let $p(\cdot)\in\mathcal{P}^{\log}(\mathbb B)$ be such that $p_->0.$  Let $w\in \Lambda$. Then
			\[\|\chi_B\|_{p(\cdot),w}\simeq w(B)^{\frac{1}{p_{+}(B)}}\simeq w(B)^{\frac{1}{p_{-}(B)}}\simeq w(B)^{\frac{1}{p(x)}}\simeq w(B)^{\frac{1}{p_{B}}}\] for all pseudo-balls $B$ of $\mathbb B$ such that $w(B)>0$ and for all $x\in B$. 
		\end{lem}
		\begin{proof}
			Take $\mathbb B$ for $B$ and $B$ for $E$ in Definition \ref{def-lambda}. We have
			$$\left (C^{-1}\frac {\mu_\alpha (B)}{\mu_\alpha (\mathbb B)} \right )^{\frac 1\delta}w(\mathbb B)\leq w(B)\leq w(\mathbb B).$$
			So
			\[w(\mathbb B)^{p-(B)-p_+(B)}\leq w(B)^{p-(B)-p_+(B)}\lesssim \mu_\alpha (B)^{\frac 1\delta(p-(B)-p_+(B))} w(\mathbb B)^{p-(B)-p_+(B)}.\]
		It is easy to check that
		\[\min (1, w(\mathbb B)^{p_--p_+})\leq w(\mathbb B)^{p-(B)-p_+(B)}\leq \max (1, w(\mathbb B)^{p_--p_+}).\]
		Next, combining with Lemma \ref{lee3} gives
		\[\min (1, w(\mathbb B)^{p_--p_+})\leq w(B)^{p-(B)-p_+(B)}\lesssim C_\delta \max (1, w(\mathbb B)^{p--p_+}).\]
		We have thus proved the estimates $w(B)^{\frac{1}{p_{+}(B)}}\simeq w(B)^{\frac{1}{p_{-}(B)}}\simeq w(B)^{\frac{1}{p(x)}}\simeq w(B)^{\frac{1}{p_{B}}}$ for all $x\in B.$
				
		On the other hand, from Lemma~\ref{le8} we have
		\begin{equation*}
		\min\left(w(B)^{\frac{1}{p_-(B)}}, w(B)^{\frac{1}{p_+(B)}}\right)
		\leq\|\chi_B\|_{p(\cdot),w}\leq\max\left(w(B)^{\frac{1}{p_-(B)}},w(B)^{\frac{1}{p_+(B)}}\right).
		\end{equation*}
		Hence \[\|\chi_B\|_{p(\cdot),w}\simeq w(B)^{\frac{1}{p_-(B)}}\simeq w(B)^{\frac{1}{p_+(B)}}.\]
	\end{proof}

		We recall again the  notation $w'(y)=w(y)^{1-p'(y)}.$
\begin{lem}\label{l0}
Let $p(\cdot)\in\mathcal{P}^{\log}_{\pm}(\mathbb B)$ and $w\in B^{+}_{p(\cdot)}.$ Then 
\begin{equation}\|w^{-1}\chi_B\|_{\frac {p'(\cdot)}{p(\cdot)}}\simeq\left(\rho_{\frac{p'(\cdot)}{p(\cdot)}}(w^{-1}\chi_B)\right)^{p_B-1}=w'(B)^{p_B-1}.\label{eq17}\end{equation}
\end{lem}

\begin{proof}
Let $w\in B^{+}_{p(\cdot)}$ and $B\in\mathcal{B}$. By definition we have 
			\begin{equation}
			\frac{1}{\mu_{\alpha}(B)^{p_B}}w(B)\|w^{-1}\chi_B\|_{\frac {p'(\cdot)}{p(\cdot)}}\leq[w]_{B^{+}_{p(\cdot)}}.\label{eq14}
			\end{equation} 	
			On the other hand, by the H\"older inequality (Proposition~\ref{hold}, assertion 1) and Lemma~\ref{le9}, we have
			\[\mu_{\alpha}(B)=\int_{B}w(y)^{\frac{1}{p(y)}}w(y)^{-\frac{1}{p(y)}}d\mu_{\alpha}(y)\leq2\|w^{\frac{1}{p(\cdot)}}\chi_B\|_{p(\cdot)}\|w^{-\frac{1}{p(\cdot)}}\chi_B\|_{p'(\cdot)}\simeq w(B)^{\frac{1}{p_B}}\|w^{-\frac{1}{p(\cdot)}}\chi_B\|_{p'(\cdot)}.\]
			Hence \begin{equation}\left \Vert\frac{w(B)^{\frac{1}{p_B}}}{\mu_{\alpha}(B)}w^{-\frac{1}{p(\cdot)}}\chi_B\right \Vert_{p'(\cdot)}\gtrsim1.\label{eq13}\end{equation}
			\\Consequently, from \eqref{eq13}, Lemma~\ref{le8}, Lemma~\ref{le9} and Corollary \ref{rqq2}, we have:
			\begin{align*}
			1\lesssim\rho_{p'(\cdot)}\left(\frac{w(B)^{\frac{1}{p_B}}}{\mu_{\alpha}(B)}w^{-\frac{1}{p(\cdot)}}\chi_B\right)&=\int_B\left(\frac{w(B)^{\frac{1}{p_B}}}{\mu_{\alpha}(B)}\right)^{p'(y)}w(y)^{-\frac{p'(y)}{p(y)}}d\mu_{\alpha}(y)\\
			&\simeq\int_B\left(\frac{w(B)}{\mu_{\alpha}(B)^{p_B}}\right)^{\frac{p'(y)}{p(y)}}w(y)^{-\frac{p'(y)}{p(y)}}d\mu_{\alpha}(y)\\
			&=\rho_{\frac{p'(\cdot)}{p(\cdot)}}\left(\frac{w(B)}{\mu_{\alpha}(B)^{p_B}}w^{-1}\chi_B\right).
			\end{align*}
			So by Lemma~\ref{le8}, we have 
			\begin{equation}
			\left \Vert\frac{w(B)}{\mu_{\alpha}(B)^{p_B}}w^{-1}\chi_B\right \Vert_{\frac {p'(\cdot)}{p(\cdot)}}\gtrsim1.\label{eq15}
			\end{equation}
			Thus from \eqref{eq14} and \eqref{eq15}, we have
			\begin{equation}\label{eq16}
			\|w^{-1}\chi_B\|_{\frac {p'(\cdot)}{p(\cdot)}}\simeq\frac{\mu_{\alpha}(B)^{p_B}}{w(B)}.
			\end{equation}
			\\Furthermore, from Remark~\ref{rq2} and as $p(\cdot)\in \mathcal P_\pm,$  we have the equivalences
			\begin{align*}\mu_{\alpha}(B)^{p_+(B)}\simeq\mu_{\alpha}(B)^{p_-(B)}&\iff\mu_{\alpha}(B)^{\frac{1}{p_-(B)-1}}\simeq\mu_{\alpha}(B)^{\frac{1}{p_+(B)-1}}\\
			&\iff\mu_{\alpha}(B)^{\frac{p_B}{p_-(B)-1}}\simeq\mu_{\alpha}(B)^{\frac{p_B}{p_+(B)-1}},
			\end{align*} 
from which we deduce that
			\[\mu_{\alpha}(B)^{\frac{p_B}{p_+(B)-1}}\simeq\mu_{\alpha}(B)^{\frac{p_B}{p_-(B)-1}}\simeq\mu_{\alpha}(B)^{\frac{p_B}{p_B-1}}\]
since $p_-(B)\leq p_B\leq p_+(B).$
			Similarly, by Lemma~\ref{le9}, we deduce from the estimate 
			\[w(B)^{p_B}\simeq w(B)^{p_+(B)}\simeq w(B)^{p_-(B)}\]
			that
			\[w(B)^{\frac{1}{p_-(B)-1}}\simeq w(B)^{\frac{1}{p_+(B)-1}}\simeq w(B)^{\frac{1}{p_B-1}}.\] 
			So from \eqref{eq16} we have 
			\[\|w^{-1}\chi_B\|_{\frac {p'(\cdot)}{p(\cdot)}}^{\frac{1}{p_+(B)-1}}\simeq\|w^{-1}\chi_B\|_{\frac {p'(\cdot)}{p(\cdot)}}^{\frac{1}{p_-(B)-1}}\simeq\|w^{-1}\chi_B\|_{\frac {p'(\cdot)}{p(\cdot)}}^{\frac{1}{p_B-1}}.\]
Since $\rho_{\frac{p'(\cdot)}{p(\cdot)}}\left(w^{-1}\chi_B\right )=w'(B),$ combining with Lemma~\ref{le8} where $\frac {p'(\cdot)}{p(\cdot)}$ replaces $p(\cdot),$  we obtain the required result.
\end{proof}

		\begin{prop}\label{p1}
			Let $p(\cdot)\in\mathcal{P}^{\log}_{\pm}(\mathbb B)$ and $w\in B^{+}_{p(\cdot)}.$ Then $w'\in B^{+}_{p'(\cdot)}$ and
			\[\left (\frac{1}{\mu_{\alpha}(B)^{p'_B}}\|w'\chi_{B}\|_1\|w'^{-1}\chi_B\|_{\frac {p(\cdot)}{p'(\cdot)}}\right )^{p_B-1}\simeq \frac{1}{\mu_{\alpha}(B)^{p_B}}\|w\chi_{B}\|_1\|w^{-1}\chi_B\|_{\frac {p'(\cdot)}{p(\cdot)}}\simeq\frac{w(B)}{\mu_{\alpha}(B)}\left(\frac{w'(B)}{\mu_{\alpha}(B)}\right)^{p_B-1}\]
			for all pseudo-balls $B\in \mathcal B.$
		\end{prop}
		\begin{proof}
			We recall that $w\in \Lambda$ by Remark ~\ref{embed}. Hence, from Lemma~\ref{le9} used with $\frac {p(\cdot)}{p'(\cdot)}$ replacing $p(\cdot),$  equation \eqref{eq17} and the property $w\in B^{+}_{p(\cdot)}$, we obtain
			\begin{align}
			\frac{1}{\mu_{\alpha}(B)^{p'_B}}w'(B)\|w'^{-1}\chi_B\|_{\frac {p(\cdot)}{p'(\cdot)}}&=\frac{1}{\mu_{\alpha}(B)^{p'_B}}w'(B)\|\chi_B\|_{\frac {p(\cdot)}{p'(\cdot)}, w}\nonumber\\ 
			&\simeq\frac{1}{\mu_{\alpha}(B)^{p'_B}}w'(B)w(B)^{\frac{1}{p_B-1}}\nonumber\\
			&=\left(\frac{w(B)}{\mu_{\alpha}(B)^{p_B}}w'(B)^{p_B-1}\right)^{\frac{1}{p_B-1}}\label{eq18}\\
			&\simeq\left(\frac{w(B)}{\mu_{\alpha}(B)^{p_B}}\|w^{-1}\chi_B\|_{\frac {p'(\cdot)}{p(\cdot)}}\right)^{\frac{1}{p_B-1}}\label{eq19}\\
			&\leq [w]_{B^{+}_{p(\cdot)}}^{\frac{1}{p_B-1}}.\nonumber
			\end{align}
			Hence $w'\in B^{+}_{p'(\cdot)}$ and from \eqref{eq18} and \eqref{eq19}, we deduce that
			\[\frac{w(B)}{\mu_{\alpha}(B)^{p_B}}\|w^{-1}\chi_B\|_{\frac {p'(\cdot)}{p(\cdot)}}\simeq\frac{w(B)}{\mu_{\alpha}(B)}\left(\frac{w'(B)}{\mu_{\alpha}(B)}\right)^{p_B-1
			}.\] 
		\end{proof}

\begin{defn}
		Let $p(\cdot)\in\mathcal{P}_{\pm}(\mathbb{B})$ and let $w$ be a weight. We say that $w$ is in the $B^{++}_{p(\cdot)}$ class if 
		\begin{equation}[w]_{B^{++}_{p(\cdot)}}:=\sup \limits_{B\in \mathcal B}\frac{w(B)}{\mu_{\alpha}(B)}\left(\frac{w'(B)}{\mu_{\alpha}(B)}\right)^{p_B-1}<\infty.\label{eq2-2}\end{equation}
	\end{defn}
It is easy to check the following Proposition. 
\begin{prop}\label{B-equiv}
	 Let $p(\cdot)\in\mathcal{P}_{\pm}(\mathbb{B}).$ The following two assertions are equivalent.
\begin{enumerate}
\item[1.] $w\in B^{++}_{p(\cdot)};$ 
\item[2.]  $w'\in B^{++}_{p'(\cdot)}$.
\end{enumerate}
\end{prop}
\begin{lem}\label{lee3-11}
	Let $p(\cdot)\in\mathcal{P}^{\log}_{\pm}(\mathbb B)$, and $w\in B_{p(\cdot)}^{++}$. Then for all pseudo-balls $B$ of $\mathbb{B}$,
	\[\|\chi_B\|_{p(\cdot),w}\simeq w(B)^{\frac{1}{p_{+}(B)}}\simeq w(B)^{\frac{1}{p_{-}(B)}}\simeq w(B)^{\frac{1}{p_B}}.\]
\end{lem}
\begin{proof}
	Since $w\in B_{p(\cdot)}^{++}$, we have $w'(\mathbb{B})<\infty$ and it follows from Lemma~\ref{le8} that $\|\chi_{\mathbb B}\|_{p'(\cdot),w'}<\infty$. Hence, by the H\"older inequality, we obtain
	\begin{align*}
	\mu_{\alpha}(B)&\leq2\|\chi_B\|_{p(\cdot),w}\|\chi_{B}\|_{p'(\cdot),w'}\nonumber\\
	&\leq2\|\chi_B\|_{p(\cdot),w}\|\chi_{\mathbb B}\|_{p'(\cdot),w'}.
	\end{align*}
	Therefore, from Lemma~\ref{lee3}, we get
	\begin{align}
	\|\chi_B\|_{p(\cdot),w}^{p_{-}(B)-p_{+}(B)}&\lesssim\mu_{\alpha}(B)^{p_{-}(B)-p_{+}(B)}\|\chi_{\mathbb B}\|_{p'(\cdot),w'}^{p_{+}(B)-p_{-}(B)}\nonumber\\
	&\lesssim\max\left (1,\|\chi_{\mathbb B}\|_{p'(\cdot),w'}^{p_{+}-p_{-}}\right ).\label{eqqn1}
	\end{align}  
	On the other hand, using again $w\in B_{p(\cdot)}^{++}$, we have $w(\mathbb{B})<\infty$ and hence $\|\chi_{\mathbb B}\|_{p(\cdot),w}<\infty$. Then
	\begin{equation}
	\|\chi_B\|_{p(\cdot),w}^{p_{+}(B)-p_{-}(B)}\leq \left \Vert\chi_{\mathbb B}\right \Vert_{p(\cdot),w}^{p_{+}(B)-p_{-}(B)}\leq\max\left (1,\|\chi_{\mathbb B}\|_{p(\cdot),w}^{p_{+}-p_{-}}\right ).\label{eqqn2}
	\end{equation} 
	Thus from \eqref{eqqn1} and \eqref{eqqn2}, we have
	\[\|\chi_B\|_{p(\cdot),w}^{p_{+}(B)}\simeq\|\chi_B\|_{p(\cdot),w}^{p_{-}(B)},\]
	and from Lemma~\ref{le8} we conclude the proof. 
\end{proof}
\begin{lem}\label{prop0}
	Let $p(\cdot)\in\mathcal{P}^{\log}_{\pm}(\mathbb B).$ Then  $B_{p(\cdot)}^{++}\subset \Lambda$.
\end{lem}
\begin{proof}
	Let $w\in B_{p(\cdot)}^{++}.$ Let $B\in \mathcal B$ and let $E$ be a measurable subset of $B.$ By the H\"older inequality and from Lemma~\ref{le8}, we have	
	\begin{align*}
	\mu_{\alpha}(E)&\leq2\|\chi_E\|_{p(\cdot),w}\|\chi_{E}\|_{p'(\cdot),w'}\nonumber\\
	&\leq2\|\chi_E\|_{p(\cdot),w}\|\chi_{B}\|_{p'(\cdot),w'}\nonumber\\
	&\leq2\max\left(w(E)^{\frac{1}{p_+(B)}},w(E)^{\frac{1}{p_-(B)}}\right)\|\chi_{B}\|_{p'(\cdot),w'}.
	\end{align*}
However, since $w'\in B_{p'(\cdot)}^{++}$ by Proposition~\ref{B-equiv},  from Lemma~\ref{lee3-11} we have
	\begin{equation*}
	\|\chi_B\|_{p'(\cdot),w'}\simeq w'(B)^{\frac{1}{p'_+(B)}}\simeq w'(B)^{\frac{1}{p'_-(B)}}\simeq w'(B)^{\frac{1}{p'_B}}.\label{eqaa1}
	\end{equation*}
	Hence using $w\in B_{p(\cdot)}^{++},$ we deduce that,
	\begin{align*}
		\mu_{\alpha}(E)&\lesssim\max\left(w(E)^{\frac{1}{p_+(B)}}w'(B)^{\frac{1}{p'_+(B)}},w(E)^{\frac{1}{p_-(B)}}w'(B)^{\frac{1}{p'_-(B)}}\right)\\
		&\leq\max\left ([w]_{B_{p(\cdot)}^{++}}^{\frac{1}{p_+(B)}},[w]_{B_{p(\cdot)}^{++}}^{\frac{1}{p_-(B)}}\right )\max\left(\left(\frac{w(E)}{w(B)}\right)^{\frac{1}{p_+(B)}},\left(\frac{w(E)}{w(B)}\right)^{\frac{1}{p_-(B)}}\right)\mu_{\alpha}(B)\\
		&\leq\max\left ([w]_{B_{p(\cdot)}^{++}}^{\frac{1}{p_+}},[w]_{B_{p(\cdot)}^{++}}^{\frac{1}{p_-}}\right )\left(\frac{w(E)}{w(B)}\right)^{\frac{1}{p_+(B)}}\mu_{\alpha}(B)\\
		&\leq\max\left ([w]_{B_{p(\cdot)}^{++}}^{\frac{1}{p_+}},[w]_{B_{p(\cdot)}^{++}}^{\frac{1}{p_-}}\right )\left(\frac{w(E)}{w(B)}\right)^{\frac{1}{p_+}}\mu_{\alpha}(B).
	\end{align*}
	Therefore, 
	\begin{equation*}
	\frac{\mu_{\alpha}(E)}{\mu_{\alpha}(B)}\lesssim\left(\frac{w(E)}{w(B)}\right)^{\frac{1}{p_+}}.
	\end{equation*} 
\end{proof}
		\begin{prop}\label{rdefi4}
Let $p(\cdot)\in \mathcal P^{\log}_\pm (\mathbb B).$ Then $B^{+}_{p(\cdot)}= B^{++}_{p(\cdot)}.$			
		\end{prop}

\begin{proof}
For the inclusion $B^{+}_{p(\cdot)}\subset B^{++}_{p(\cdot)},$ apply Proposition~\ref{p1}. For the reverse inclusion, apply Lemma~\ref{prop0} and Lemma~\ref{le9}. 
\end{proof}

		\begin{rem}\label{rqm11}
			From Proposition~\ref{p1}, Remark~\ref{embed} and Lemma~\ref{le9}, we have the inclusion $B^{+}_{p(\cdot)}\subset B_{p(\cdot)}$ for $p(\cdot)\in \mathcal P^{\log}_\pm (\mathbb B).$
		\end{rem}
		
		Now we prove the reverse inclusion. In this direction, we first state the following result. 
		\begin{lem}\label{lee3-1}
			Let $p(\cdot)\in\mathcal{P}^{\log}_{\pm}(\mathbb{B})$ and $w\in B_{p(\cdot)}$. Then there exists a constant $C>1$
			\[\frac 1C\le\|\chi_B\|_{p(\cdot),w}^{p_{-}(B)-p_{+}(B)}\leq C\]
for all $B\in \mathcal B.$
		\end{lem}
		\begin{proof}
			By the H\"older inequality, we have
			\begin{equation*}
			\mu_{\alpha}(B)\leq2\|\chi_B\|_{p(\cdot),w}\|\chi_{B}\|_{p'(\cdot),w'}\label{eqq1}
			\end{equation*}
			 and as $w\in B_{p(\cdot)}$, from Lemma~\ref{lee3} and according to the estimate $\|\chi_{\mathbb B}\|_{p'(\cdot),w'}<\infty$ given by Remark~\ref{elem}, we have
			 \begin{align*}
			 \|\chi_B\|_{p(\cdot),w}^{p_{-}(B)-p_{+}(B)}&\lesssim\mu_{\alpha}(B)^{p_{-}(B)-p_{+}(B)}\|\chi_{B}\|_{p'(\cdot),w'}^{p_{+}(B)-p_{-}(B)}\\
			 &\simeq\|\chi_{\mathbb B}\|_{p'(\cdot),w'}^{p_{+}(B)-p_{-}(B)}\\
			 &\lesssim\max\left (1,\|\chi_{\mathbb{B}}\|_{p'(\cdot),w'}^{p_{+}-p_{-}}\right ).
			 \end{align*}
			 
			On the other hand, according to the estimate $\|\chi_{\mathbb B}\|_{p(\cdot),w}<\infty$ given by Remark~\ref{elem}, we have 
			\begin{equation*}
			\|\chi_B\|_{p(\cdot),w}^{p_{+}(B)-p_{-}(B)}\leq\max(1,\|\chi_{\mathbb{B}}\|_{p(\cdot),w}^{p_{+}-p_{-}}).
			\end{equation*} 
		\end{proof}
		From Lemma~\ref{le8} and Lemma~\ref{lee3-1}, we deduce the following corollary.
		\begin{cor}\label{core1}
			Let $p(\cdot)\in\mathcal{P}^{\log}_{\pm}(\mathbb{B})$, and $w\in B_{p(\cdot)}$. Then 
			\[\|\chi_B\|_{p(\cdot),w}\simeq w(B)^{\frac{1}{p_{+}(B)}}\simeq w(B)^{\frac{1}{p_{-}(B)}}\simeq w(B)^{\frac{1}{p_{B}}}.\]
for all $B\in\mathcal{B}$,
		\end{cor}

		We next state the following theorem.
		\begin{thm}\label{equiv-def}
			Let $p(\cdot)\in\mathcal{P}^{log}_{\pm}(\mathbb B).$  Then $B_{p(\cdot)}=B^{+}_{p(\cdot)}=B^{++}_{p(\cdot)}$.
		\end{thm}
		\begin{proof}
			From Remark~\ref{rqm11}, we have $B^+_{p(\cdot)}\subset B_{p(\cdot)}$.
			Let $w\in B_{p(\cdot)}.$ By Proposition~\ref{w-w'}, $w'\in B_{p'(\cdot)}.$ It follows from Corollary~\ref{core1} that
			\begin{equation*}
			\frac{w(B)}{\mu_{\alpha}(B)}\left(\frac{w'(B)}{\mu_{\alpha}(B)}\right)^{p_B-1}\simeq\left(\frac{1}{\mu_{\alpha}(B)}\|\chi_B\|_{p(\cdot),w}\|\chi_B\|_{p'(\cdot),w'}\right)^{p_B}\leq\max \left (1,[w]_{B_{p(\cdot)}}^{p_+}\right ).
			\end{equation*}
		\end{proof}
		To end this section, we record with the same proof the following analogous theorem for the variable Muckenhoupt weight classes.
		\begin{thm}\label{th-Ap}
			Let $w$ be a weight and let $p(\cdot)\in \mathcal{P}^{\log}_{\pm}(\mathbb{B})$. The following three assertions are equivalent.
			\begin{enumerate}
				\item[1.] $w\in A_{p(\cdot)}$;
				\item[2.] 
				$\sup \limits_{B}\frac 1{\mu_{\alpha}(B)^{p_B}}\|w\chi_{B}\|_1\|w^{-1}\chi_B\|_{\frac {p'(\cdot)}{p(\cdot)}}<\infty,$ where the $\sup$ is taken over all pseudo-balls of $\mathbb B;$
				\item[3.] 
				$\sup \limits_{B}\frac{w(B)}{\mu_{\alpha}(B)}\left(\frac{w'(B)}{\mu_{\alpha}(B)}\right)^{p_B-1}<\infty,$ where the $\sup$ is taken over all pseudo-balls of $\mathbb B.$
			\end{enumerate}
		\end{thm}
		\section{Proof of the necessary condition in Theorem~\ref{th1}}\label{sec-6}
The aim of this section is to prove the following result.
	\begin{prop}\label{p2}
		Let $w$ be a weight and let $p(\cdot)\in\mathcal{P}^{\log}_{\pm}(\mathbb B).$ If the Bergman projector is bounded on $L^{p(\cdot)}(wd\mu_{\alpha})$, then $w\in B_{p(\cdot)}$.
	\end{prop}
	\begin{proof}
According to Theorem~\ref{equiv-def}, it suffices to prove that $w\in B_{p(\cdot)}^{++},$ i.e.
the following estimate holds
\begin{equation}\label{nec}
\sup \limits_{B\in \mathcal B}\frac{w(B)}{\mu_{\alpha}(B)}\left(\frac{w'(B)}{\mu_{\alpha}(B)}\right)^{p_B-1}<\infty.
\end{equation}
	From Proposition~\ref{l1}, we have that $w^{-\frac{1}{p(\cdot)}}\in L^{p'(\cdot)}(d\mu_\alpha)$ and  from Proposition~\ref{l2} we have $w^{\frac{1}{p(\cdot)}}\in L^{p(\cdot)}(d\mu_\alpha).$ In particular, $w(\mathbb B)<\infty$ and $w'(\mathbb B)<\infty.$
		Thus we just have to show the estimate \eqref{nec} for the pseudo-balls of radius smaller than a positive constant $R_0$, because if the radius of $B$ is  larger, then $B$ can be identified with $\mathbb B$. 
We shall use the following lemma.
	\begin{lem}\cite{ref3}\label{le7}
There exist three positive numbers $R_0, c$ and $C_\alpha$ such that the following holds.		For every pseudo-ball $B^1\in \mathcal B$  of  radius $R<R_0,$ there exists a pseudo-ball $B^2\in \mathcal B$ of same radius such that $d(B^1, B^2)=cR$, that satisfies the following property:  for every non-negative measurable function $f$  supported in $B^i$ and for two distinct superscripts $i,j\in \{1,2\},$ we have
		\begin{equation}|
  P_{\alpha} f|\geq C_{\alpha}\chi_{B^j}\mu_{\alpha}(B^i)^{-1}\int_{B^i}fd\mu_{\alpha}\label{eq3}.\end{equation}
	\end{lem}
	Thus, by taking $f=\chi_{B^i}$ in \eqref{eq3} we obtain :
	\begin{equation*}
	|P_{\alpha}\chi_{B^i}(z)|\geq\chi_{B^j}(z) C_{\alpha}\mu_{\alpha}(B^i)^{-1}\int_{B^i}\chi_{B^i}d\mu_{\alpha}\simeq \chi_{B^j}(z).
	\end{equation*}
	Using the  growth of the norm $\|\cdot\|_{p(\cdot), w},$ we obtain
	\begin{equation*}
	\|P_{\alpha}\chi_{B^i}\|_{p(\cdot),w}=\|w^{\frac{1}{p(\cdot)}}P_{\alpha}\chi_{B^i}\|_{p(\cdot)}\gtrsim \|\chi_{B^j}w^{\frac{1}{p(\cdot)}}\|_{p(\cdot)}
	\end{equation*}
	So using the fact that $P_{\alpha}$ is bounded on $L^{p(\cdot)}(wd\mu_{\alpha})$,  we obtain:
	\begin{equation*}
	\|\chi_{B^j}w^{\frac{1}{p(\cdot)}}\|_{p(\cdot)}\lesssim\|P_{\alpha}\|\|\chi_{B^i}\|_{p(\cdot),w}
	\end{equation*}
	We then deduce that 
	\begin{equation}
	\|\chi_{B^1}w^{\frac{1}{p(\cdot)}}\|_{p(\cdot)}\simeq \|\chi_{B^2}w^{\frac{1}{p(\cdot)}}\|_{p(\cdot)}.\label{eq7}
	\end{equation}
	\\In the rest of the proof, we shall take $f=w'\chi_{B^1}.$ We have $f\in L^{p(\cdot)}(wd\mu_{\alpha})$ since
\[\rho_{p(\cdot), w}(f)=\int_{\mathbb B}w(z)^{-p'(z)}\chi_{B^1}(z)w(z)d\mu_{\alpha}(z)= \rho_{p'(\cdot)}(w^{-\frac{1}{p(\cdot)}}\chi_{B^1})<\infty\]
by Proposition \ref{l1}. However, $\rho_{p(\cdot), w}(f)=\int_{B^1} w'd\mu_\alpha=w'(B^1).$
Also, from \eqref{eq3} and the previous equality,  
we have
	\[\chi_{B^2}(z)w'(B^1)\leq C_{\alpha}^{-1} \mu_{\alpha}(B^1) |P_{\alpha}f(z)|.\]
	Moving to the norm $\|\cdot\|_{p(\cdot), w},$ we obtain 
	\begin{equation*}
	\|w^{\frac{1}{p(\cdot)}}\chi_{B^2}\|_{p(\cdot)}w'(B^1)\leq C_{\alpha}^{-1}\mu_{\alpha}(B^1)\|P_{\alpha}f\|_{p(\cdot),w}
	\end{equation*}
	Then using the boundedness of $P_{\alpha}$ on $L^{p(\cdot)}(wd\mu_{\alpha}),$ the previous inequality implies  
	\begin{equation*}
	\|w^{\frac{1}{p(\cdot)}}\chi_{B^2}\|_{p(\cdot)}w'(B^1)\leq C_{\alpha}^{-1}\mu_{\alpha}(B^1)\|P_{\alpha}\|\|f\|_{p(\cdot),w}
	\end{equation*}
	and combining with \eqref{eq7}, we obtain the following lemma.
\begin{lem}[Main Lemma]\label{lem}
Suppose that $P_\alpha$ is bounded on $L^{p(\cdot)}(wd\mu_{\alpha}).$ Then
	\begin{equation}
	\|w^{\frac{1}{p(\cdot)}}\chi_B\|_{p(\cdot)}w'(B)\leq CC_{\alpha}^{-1}\mu_{\alpha}(B)\|P_{\alpha}\|\|w'\chi_{B}\|_{p(\cdot),w}	\label{eq22}
	\end{equation}
for every pseudo-ball $B\in \mathcal B$ of radius smaller than $R_0.$ The absolute constants $R_0, C$ and $C_\alpha$ were respectively defined in Lemma~\ref{le7}, \eqref{eq7} and \eqref{eq3}.
\end{lem}
	At this level, we need to calculate $\|w^{\frac{1}{p(\cdot)}}\chi_B\|_{p(\cdot)}$ and $\|\|w'\chi_{B}\|_{p(\cdot),w}.$  This calculation is not as obvious as in the case where $p(\cdot)$ is constant. 
	\begin{lem}\label{P-P'}
		Let $p(\cdot)\in\mathcal{P}(\mathbb B)$. If $P_\alpha$ is bounded on $L^{p(\cdot)}(w)$, then $P_\alpha$ is bounded on $L^{p'(\cdot)}(w')$.
	\end{lem}
\begin{proof}
We first recall that the weighted Bergman projector $P_\alpha$ is the orthogonal projector from the (Hilbert-) Lebesgue space $L^2 (d\mu_\alpha)$ to its closed subspace $L^2 (d\mu_\alpha)\cap Hol (\mathbb B)$ (the standard weighted Bergman space). We call $\mathcal C_c (\mathbb B)$ the space of continuous functions  with compact support in $\mathbb B.$ By Proposition \ref{dense}, $\mathcal C_c (\mathbb B)$ is a dense subspace of $L^{p(\cdot)}(w)$ and $L^{p'(\cdot)}(w').$	From Proposition~\ref{egalm} and the boundedness of $P_\alpha$  on $L^{p(\cdot)}(w),$ for all  $f\in \mathcal C_c (\mathbb B),$  we have
	\begin{align*}
		\|P_\alpha f\|_{p'(\cdot),w'}&=\sup_{g\in \mathcal C_c (\mathbb B:\|g\|_{p(\cdot),w}=1}\left|\int_{\mathbb B}P_\alpha f(\zeta)\overline{g(\zeta)}d\mu_{\alpha}(\zeta)\right|\\
		&=\sup_{g\in \mathcal C_c (\mathbb B): \|g\|_{p(\cdot),w}=1}\left|\int_{\mathbb B}f(\zeta) \overline{P_\alpha g(\zeta)}d\mu_{\alpha}(\zeta)\right|\\
		&\leq 2\sup_{g\in \mathcal C_c (\mathbb B): \|g\|_{p(\cdot),w}=1}\|f\|_{p'(\cdot),w'}\|P_\alpha g\|_{p(\cdot),w}\\
		&\leq 2\left \Vert P_\alpha \right \Vert\|f\|_{p'(\cdot),w'}.
	\end{align*}
We have used the elementary fact that $\mathcal C_c (\mathbb B)$ is contained in $L^2 (d\mu_\alpha).$ For the last but one inequality, we used the H\"older inequality.
\end{proof}
\begin{lem}\label{weak_cv}
	Let $p(\cdot)\in\mathcal{P}(\mathbb B)$ and let $w$ be a weight. If $P_\alpha$ is bounded on $L^{p(\cdot)}(w)$, then for all $t>0$,
	\[\|t\chi_{\{|P_\alpha f|> t\}}\|_{p(\cdot),w}\leq \left \Vert P_\alpha \right \Vert\|f\|_{p(\cdot),w}.\]
\end{lem}
\begin{proof}
	It suffices to remark that for all $t>0,\;t\chi_{\{|P_\alpha f|> t\}}\leq |P_\alpha f|.$
\end{proof}

\begin{lem}\label{leP-in}
	Let $p(\cdot)\in\mathcal{P}^{\log}_{\pm}(\mathbb B)$ and $w$ be a weight. If $P_\alpha$ is bounded on $L^{p(\cdot)}(w)$, then
	\[\|\chi_B\|_{p(\cdot),w}\simeq w(B)^{\frac{1}{p_{+}(B)}}\simeq w(B)^{\frac{1}{p_{-}(B)}}\] for all pseudo balls $B$ of $\mathbb B$.
\end{lem}
\begin{proof}
	If $\|\chi_B\|_{p(\cdot),w}\geq1,$ then $\left \Vert \chi_{\mathbb B}\right \Vert_{p(\cdot),w}^{p_--p_{+}}\leq \|\chi_B\|_{p(\cdot),w}^{p_--p_{+}}\leq \|\chi_B\|_{p(\cdot),w}^{p_{-}(B)-p_{+}(B)}\leq 1.$ So
\begin{equation}\label{simeq}
\|\chi_B\|_{p(\cdot),w}^{p_{-}(B)-p_{+}(B)}\simeq1. 
\end{equation}
 Otherwise, if $\|\chi_B\|_{p(\cdot),w}<1,$ then by the H\"older inequality,  we have
	\begin{align}
	\mu_{\alpha}(B)&\leq2\|\chi_B\|_{p(\cdot),w}\|\chi_{B}\|_{p'(\cdot),w'}\nonumber\\
	&\leq2\|\chi_B\|_{p(\cdot),w}\|\chi_{\mathbb B}\|_{p'(\cdot),w'}\label{eqqq1}
	\end{align}
	Hence from \eqref{eqqq1} and Lemma~\ref{lee3}, we have 
	\begin{align}
	\|\chi_B\|_{p(\cdot),w}^{p_{-}(B)-p_{+}(B)}&\leq 2^{p_+-p_-}\mu_{\alpha}(B)^{p_{-}(B)-p_{+}(B)}\|\chi_{\mathbb B}\|_{p'(\cdot),w'}^{p_{+}(B)-p_{-}(B)}\nonumber\\
	&\lesssim \max \left (1, \|\chi_{\mathbb B}\|_{p'(\cdot),w'}^{p_{+}-p_{-}}\right )\label{eqqq3}.
	\end{align} 
We point out that  $\|\chi_{\mathbb B}\|_{p'(\cdot),w'}<\infty$ according to  Proposition~\ref{l1}, since $P_\alpha$ is bounded on $L^{p(\cdot)}(w).$ 
	On the other hand, using again the boundedness of $P_\alpha$ on $L^{p(\cdot)}(w),$ we have the estimate $\|\chi_{\mathbb B}\|_{p(\cdot), w}<\infty$ according to Proposition~\ref{l2}. Then
	\begin{equation}
	\|\chi_B\|_{p(\cdot),w}^{p_{+}(B)-p_{-}(B)}\lesssim \max \left (1,\|\chi_{\mathbb B}\|_{p(\cdot),w}^{p_{+}-p_{-}}\right )<\infty.\label{eqqq4}
	\end{equation}
	Thus, from \eqref{simeq}, \eqref{eqqq3} and \eqref{eqqq4}, we deduce that
	\begin{equation}
		\|\chi_B\|_{p(\cdot),w}^{p_{+}(B)}\simeq \|\chi_B\|_{p(\cdot),w}^{p_{-}(B)}
	\end{equation} 
	for all pseudo-balls of $\mathbb B.$ Applying Lemma~\ref{le8} gives
	\begin{equation*}
	\|\chi_B\|_{p(\cdot),w}\simeq w(B)^{\frac{1}{p_+(B)}}\simeq w(B)^{\frac{1}{p_-(B)}}.
	\end{equation*}
\end{proof}
\noindent
\textit{End of the proof of Proposition~\ref{p2}.}
We go back to the Main Lemma (Lemma~\ref{lem}).	On the one hand, since $P_\alpha$ is bounded on $L^{p(\cdot)}(w)$, it follows from Lemma~\ref{P-P'} that $P_\alpha$ is also bounded on $L^{p'(\cdot)}(w')$. So from Lemma~\ref {leP-in} with $p'(\cdot)$ in the place of $p(\cdot)$ and $w'$ in the place of $w,$ we have
$$\left \Vert \chi_B\right \Vert_{p'(\cdot), w'}\simeq w'(B)^{1-\frac 1{p_-(B)}}\simeq w'(B)^{1-\frac 1{p_+(B)}}.$$
This implies the estimate $w'(B)^{\frac 1{p_-(B)}}\simeq w'(B)^{\frac 1{p_+(B)}}.$ It then follows from Lemma~\ref{lem} that 
\[\|w'\chi_B\|_{p(\cdot),w}\simeq w'(B)^{\frac{1}{p_{B}}}.\]
	On the other hand,  $\|w^{\frac{1}{p(\cdot)}}\chi_{B}\|_{p(\cdot)}=\|\chi_B\|_{p(\cdot), w}\simeq w(B)^{\frac{1}{p_{B}}}$ by Lemma~\ref{leP-in}.
 The inequality \eqref{eq22} of the Main Lemma  takes the following form 
	\[w(B)^{\frac{1}{p_{B}}}w'(B)\lesssim \mu_{\alpha}(B)w'(B)^{\frac{1}{p_{B}}}.\]
Equivalently,
	\begin{equation*}
	\sup \limits_{B\in \mathcal B}\frac{w(B)}{\mu_{\alpha}(B)}\left(\frac{w'(B)}{\mu_{\alpha}(B)}\right)^{p_{B}-1}<\infty.	
	\end{equation*}
We have shown the estimate \eqref{nec}. This finishes the proof of Proposition \ref{p2}.
\end{proof}		

		\section{Boundedness on $L^{p(\cdot)} (w)$ of the maximal function $m_\alpha $ }\label{sec-7}
  
  In this section, we prove the boundedness of the maximal function $ m_\alpha $ on $L^{p(\cdot)} (w)$ when $w \in B_{p(\cdot)}$. As in \cite{ref3}, we will use the regularisation operator that we recall here with some of its properties. 
		\begin{defn}
			For all $k\in(0,1),$ we define the regularisation operator $R_k^\alpha $ of order $k$    by
			\[R_k^\alpha f(z)=\frac{1}{\mu_{\alpha}(B^k(z))}\int_{B^k(z)}f(\zeta)d\mu_{\alpha}(\zeta),\]
where $B^k(z)=\{\zeta\in\mathbb{B}:d(z,\zeta)<k(1-|z|)\}.$ 
		\end{defn}
		\begin{prop}\label{pro1} For all $k\in(0,1),$ there exists a constant $C_k>1$ such that for every non-negative locally integrable function $f,$ the following two estimates hold.
			\begin{enumerate}
				\item[1)] 
				$m_\alpha f\leq C_km_\alpha R_k^\alpha f;$
				\item[2)] 
				$C_k^{-1} m_\alpha g\leq R_k^\alpha m_\alpha g\leq C_k m_\alpha g.$
			\end{enumerate}
		\end{prop}
		\begin{lem}\label{lem1}
			Let $k\in(0,\frac{1}{2})$. If $z'\in B^k(z)$ then $z\in B^{k'}(z')$ where $k'=\frac{k}{1-k},$ and  $\chi_{B^{k}(z)}(z')\leq\chi_{B^{k'}(z')}(z)$. Moreover there exists a constant $C_k>1$ such that \[C_k^{-1}\mu_{\alpha}(B^k(z))\leq\mu_{\alpha}(B^{k'}(z'))\leq C_k\mu_{\alpha}(B^k(z)).\]  
		\end{lem}
		
		\begin{lem}\label{le11}
			Let $k\in \left (0,\frac{1}{5}\right )$ and $p(\cdot)\in\mathcal{P}^{\log}_{\pm}(\mathbb B)$. For $w\in B_{p(\cdot)}$ there exists a constant $C_k>1$ such that for all $z,z'\in\mathbb B$ such that $z'\in B^k(z)$ we have:
			\[C_k^{-1}w(B^{k}(z))\leq w(B^{k'}(z'))\leq C_k w(B^{k}(z))\]
		\end{lem}
		\begin{proof}
			We have $B^k(z)\subset B^{2k'}(z')$ and $B^{k'}(z')\subset B^{6k}(z)$. From Theorem \ref{equiv-def} and Remark \ref{embed}, we have $w\in \Lambda.$ Apply Remark \ref{doubling} to conclude. 
		\end{proof}
		\begin{lem}\label{lem2}
			Let $k\in(0,\frac{1}{2}).$ There exists a positive constant $C_k$ such that for all non-negative locally integrable $f,g,$ we have
			\[\int_{\mathbb{B}}f(\zeta)R_kg(\zeta)d\mu_{\alpha}(\zeta)\leq C_k\int_{\mathbb{B}}g(z)R_kf(z)d\mu_{\alpha}(z).\]
		\end{lem}
		We also recall the following elementary lemma.
\begin{lem}\label{lem10}
For $z\in B$ and $\zeta\in B^k (z),$ we have $\zeta \in B'.$
\end{lem}
		
		In the rest of this section, to simplify the notation, we write $\sigma=R_k^\alpha w.$
		The following result is a generalisation to the variable exponent of the analogous result in \cite[Lemma 10]{ref3}.
		
		\begin{prop}\label{lem3}
			Let $p(\cdot)\in\mathcal{P}^{log}_{\pm}(\mathbb B),\,k\in(0,\frac{1}{2})$ and $w\in B_{p(\cdot)}$. Then $R_k^\alpha w\in A_{p(\cdot)}$ with $[R_k^\alpha w]_{A_{p(\cdot)}}\lesssim [w]_{B_{p(\cdot)}}.$
		\end{prop}		
				
		\begin{proof}
			From Theorem~\ref{th-Ap}, it suffices to show that
			$$\frac{\sigma(B)}{\mu_{\alpha}(B)}\left(\frac{\sigma'(B)}{\mu_{\alpha}(B)}\right)^{p(z_0)-1}\lesssim [w]_{B_{p(\cdot)}}$$
			for every pseudo-ball $B$ of $\mathbb B.$
			
			We write $a=2k+1.$ Let  $B=B(z_0,r)$ be a pseudo-ball in $\mathbb B.$ We set $B'=B(z_0, ar).$ We distinguish two cases: 1. $B\in \mathcal B;$ 2. $B$ is not a member of $\mathcal B.$
			\begin{enumerate}
				\item[1.]
				Suppose first that  $B\in\mathcal{B}.$ We claim that there exists a positive absolute constant $C_k$  such that 
				\begin{equation}\label{sigma}
				\frac{\sigma(B)}{\mu_{\alpha}(B)}\leq C_k\frac{w(B')}{\mu_{\alpha}(B')}.
				\end{equation}
				Indeed, from the Fubini-Tonelli theorem and Lemma~\ref{lem1} we have
				\begin{align*}
				\sigma(B)&=\int_B\sigma(z)d\mu_{\alpha}(z)\\
				&=\int_B\left(\frac{1}{\mu_{\alpha}(B^k(z))}\int_{B^k(z)}w(\zeta)d\mu_{\alpha}(\zeta)\right)d\mu_{\alpha}(z)\\
				&=\int_{\mathbb B}\left(\int_{\mathbb B}\frac{\chi_{B^k(z)}(\zeta)\chi_B(z)}{\mu_{\alpha}(B^k(z))}d\mu_{\alpha}(z)\right)w(\zeta)d\mu_{\alpha}(\zeta)\\
				&\lesssim\int_{\mathbb B}\left(\int_{\mathbb B}\frac{\chi_{B^{k'}(\zeta)}(z)\chi_{B'}(\zeta)}{\mu_{\alpha}(B^{k'}(\zeta))}d\mu_{\alpha}(z)\right)w(\zeta)d\mu_{\alpha}(\zeta)\\
				&=w(B').
				\end{align*}
				For the latter inequality, we used Lemma \ref{lem10}.
				Moreover, since $B\subset B'$ and $\mu_{\alpha}(B)\simeq\mu_{\alpha}(B'),$  we obtain that 
				\[\frac{\sigma(B)}{\mu_{\alpha}(B)}\lesssim\frac{\sigma(B')}{\mu_{\alpha}(B')}.\]
				Furthermore, from the H\"older inequality and Lemma~\ref{le9} we have
				\begin{align*}
				\sigma^{-1}(z)&=\frac{\mu_{\alpha}(B^k(z))}{w(B^k(z))}\\
				&\leq\frac{2}{w(B^k(z))}\|w^{\frac{1}{p(\cdot)}}\chi_{B^k(z)}\|_{p(\cdot)}\|w^{-\frac{1}{p(\cdot)}}\chi_{B^k(z)}\|_{p'(\cdot)}\\
				&\simeq\frac{1}{w(B^k(z))}w(B^k(z))^{\frac{1}{p(z)}}w'(B^k(z))^{\frac{1}{p'(z)}}\\
				&=\left (\frac{w'(B^k(z))}{w(B^k(z))}\right )^{\frac{1}{p'(z)}}.
				\end{align*} 
				Hence
				\[\sigma'(z)=\left (\sigma^{-1}\right )^{p'(z)-1}(z)\lesssim\left(\frac{w'(B^k(z))}{w(B^k(z))}\right)^{\frac{1}{p(z)}}.\]
				From the H\"older inequality and Lemma~\ref{le9}, we have 
				\begin{align}
				\sigma'(B)&=\int_B\sigma'(z)d\mu_{\alpha}(z)\nonumber\\
				&\lesssim\int_B\left(\frac{w'(B^k(z))}{w(B^k(z))}w(z)\right)^{\frac{1}{p(z)}}w(z)^{-\frac{1}{p(z)}}d\mu_{\alpha}(z)\nonumber\\
				&\leq2\|w^{-\frac{1}{p(\cdot)}}\chi_{B}\|_{p'(\cdot)}\|\left(\frac{w'(B^k(.))}{w(B^k(.))}w(.)\right)^{\frac{1}{p(.)}}\chi_{B}\|_{p(\cdot)}\nonumber\\
				&\lesssim w'(B)^{\frac{1}{p'(z_0)}}\|\left(\frac{w'(B^k(.))}{w(B^k(.))}w(.)\right)^{\frac{1}{p(.)}}\chi_{B}\|_{p(\cdot)}\label{star}
				\end{align}
		Since $w'\in \Lambda$, from Lemma~\ref{le9} we have $\beta:=\|w'^{\frac{1}{p(\cdot)}}\chi_{B'}\|_{p(\cdot)}\simeq w'(B')^{\frac{1}{p(z_0)}}$. Thus
				as $B\subset B'$, from Lemma~\ref{lem1} and Lemma~\ref{le11} we have  
				\begin{align*}
				&\rho_{p(.)}\left(\frac{1}{\beta}\left(\frac{w'(B^k(.))}{w(B^k(.))}w\right)^{\frac{1}{p(.)}}\chi_{B}\right)\\
				&=\int_{\mathbb B}\frac{1}{\beta^{p(z)}}\frac{w'(B^k(z))}{w(B^k(z))}w(z)\chi_{B}(z)d\mu_{\alpha}(z)\\
				&\simeq\int_{\mathbb B}\frac{1}{w'(B')}\frac{w'(B^k(z))}{w(B^k(z))}w(z)\chi_{B}(z)d\mu_{\alpha}(z)\\
				&=w'(B')^{-1}
			\times\int_{\mathbb B}\left(\frac{1}{w(B^k(z))}\int_{\mathbb B} w'(\zeta)\chi_{B_{k}(z)}(\zeta)\chi_B(z)w(z)d\mu_{\alpha}(\zeta)\right)d\mu_{\alpha}(z)\\
				&\leq C_kw'(B')^{-1}
			\times\int_{\mathbb B}\left(\frac{1}{w(B_{k'}(\zeta))}\int_{\mathbb B} \chi_{B_{k'}(\zeta)}(z)\chi_{B'}(\zeta)w(z)d\mu_{\alpha}(z)\right)w'(\zeta)d\mu_{\alpha}(\zeta)\\
				&=C_k.
				\end{align*}
For the latter inequality, we used Lemma~\ref{lem10}. 
				Hence we obtain
				\[\left \Vert\left(\frac{w'(B^k(.))}{w(B^k(.))}w\right)^{\frac{1}{p(.)}}\chi_{B}\right \Vert_{p(\cdot)}\lesssim C_k w'(B')^{\frac{1}{p(z_0)}}.\]
				Consequently, we deduce from \eqref{star} that 
				\[\sigma'(B)\lesssim C_k w'(B)^{\frac{1}{p'(z_0)}}w'(B')^{\frac{1}{p(z_0)}}\leq C_k w'(B')\]because $B\subset B'.$ Moreover, as $\mu_{\alpha}(B)\simeq\mu_{\alpha}(B'),$ we have
				\[
				\frac{\sigma'(B)}{\mu_{\alpha}(B)}\leq C_k\frac{w'(B')}{\mu_{\alpha}(B')}
				\]
				and hence
				\begin{equation}\label{sigma1}
				\left(\frac{\sigma'(B)}{\mu_{\alpha}(B)}\right)^{p(z_0)-1}\leq C_k' \left(\frac{w'(B')}{\mu_{\alpha}(B')}\right)^{p(z_0)-1}.
				\end{equation}
				Combining \eqref{sigma} and \eqref{sigma1} gives
				\begin{equation}\label{Rk1}
				\frac{\sigma(B)}{\mu_{\alpha}(B)}\left(\frac{\sigma'(B)}{\mu_{\alpha}(B)}\right)^{p(z_0)-1}\leq \gamma_k\frac{w(B')}{\mu_{\alpha}(B')}\left(\frac{w'(B')}{\mu_{\alpha}(B')}\right)^{p(z_0)-1}\leq\gamma_k[w]_{B_{p(\cdot)}}
				\end{equation}
by Theorem \ref{equiv-def}.
				\item[2.]
				Suppose next that the pseudo-ball $B$ is not a member of $\mathcal B$, i.e. $r\leq1-|z_0|$. In the case where $k(1-|z_0|)\leq r\leq 1-|z_0|$, we have $B\subset B(z_0, 1-|z_0|)$ and $\mu_{\alpha}(B)\simeq(1-|z_0|)^{n+\alpha}\simeq\mu_{\alpha}(B(z_0, 1-|z_0|))$. The pseudo-ball $B(z_0, 1-|z_0|)$ is a member of $\mathcal B;$ so we can apply to it the computations of the first case. We obtain:
				
				\begin{eqnarray*}
				\frac{\sigma(B)}{\mu_{\alpha}(B)}\left(\frac{\sigma'(B)}{\mu_{\alpha}(B)}\right)^{p(z_0)-1}&\lesssim &\frac{\sigma(B(z_0, 1-|z_0|))}{\mu_{\alpha}(B(z_0, 1-|z_0|))}\left(\frac{\sigma'(B(z_0, 1-|z_0|))}{\mu_{\alpha}(B(z_0, 1-|z_0|))}\right)^{p(z_0)-1}\\
				&\lesssim & [w]_{B_{p(\cdot)}}.
				\end{eqnarray*}

				Next, if $0<r<k(1-|z_0|),$ then for $z\in B$ we have
$(1-k)(1-|z_0|)\leq1-|z|\leq(1+k)(1-|z_0|)$. This shows that $\mu_{\alpha}(B^k(z_0))\simeq\mu_{\alpha}(B^k(z)).$\\
 We also claim that $w(B^k(z_0))\simeq w(B^k(z)).$
				Indeed, it is easy to show the inclusions $B^k (z_0)\subset B(z, 4k(1-|z|))$ and $B^k (z)\subset B(z_0, 2k(2+k)(1-|z_0|)).$ The claim then follows an application of Remark \ref{doubling}. Combining with the estimate $\mu_{\alpha}(B^k(z_0))\simeq\mu_{\alpha}(B^k(z))$ gives 				
				\begin{equation}\label{sigmasim}			
				\sigma (z)\simeq \sigma (z_0)
				\end{equation}
				for every $z\in B.$
				Now, by Remark \ref{p-p'}, $p'(\cdot)$ is a member of $\mathcal P^{log}_\pm (\mathbb B).$ Then by Lemma~\ref{le9} and Corollary~\ref{rqq2}, we have $$\sigma (z_0)^{1-p'(z_0)}\simeq \sigma (z)^{1-p'(z_0)}\simeq \sigma(z)^{1-p'(z)}.$$ Combining with \eqref{sigmasim} gives 
				\begin{equation}\label{Rk2}
				\frac{\sigma(B)}{\mu_{\alpha}(B)}\left(\frac{\sigma'(B)}{\mu_{\alpha}(B)}\right)^{p(z_0)-1}\simeq 1.
				\end{equation}
				The conclusion of the lemma follows a combination of \eqref{Rk1} and \eqref{Rk2} with Theorem~\ref{equiv-def}.
			\end{enumerate}
		\end{proof}
		\begin{lem}\label{rq3}
			Let $p(\cdot)\in\mathcal{P}^{\log}_{\pm}(\mathbb B), k\in(0,\frac{1}{2})$ and $w\in B_{p(\cdot)}.$  Then 
			$$(R_k^\alpha g(z))^{p(z)}\lesssim R_k^\alpha(g^{p(\cdot)})(z)+1$$
			for all non-negative functions $g$ such that $\|g\|_{p(\cdot),w}=1$ and all $z\in \mathbb B.$
		\end{lem}
		\begin{proof}
			As $\|g\|_{p(\cdot),w}=1$, from the H\"older inequality (assertion 1 of Proposition~\ref{hold}), we have
			\begin{align*}
			\frac{1}{2\|\chi_{B^k(z)}\|_{p'(\cdot),w'}}\int_{B^k(z)}g(\zeta)d\mu_{\alpha}(\zeta)\leq\|g\chi_{B^k(z)}\|_{p(\cdot),w}=1.
			\end{align*}
			Therefore from the usual H\"older inequality  and Lemma~\ref{lee3}, we obtain 
			\begin{align*}
			(R_k^\alpha g(z))^{p(z)}&=\left(\frac{1}{2\|\chi_{B^k(z)}\|_{p'(\cdot),w'}}\int_{B^k(z)}g(\zeta)d\mu_{\alpha}(\zeta)\right)^{p(z)}\mu_{\alpha}(B^k(z))^{-p(z)}2^{p(z)}\|\chi_{B^k(z)}\|_{p'(\cdot),w'}^{p(z)}\\
			&\leq2^{p(z)}\left(\frac{1}{2\|\chi_{B^k(z)}\|_{p'(\cdot),w'}}\int_{B^k(z)}g(\zeta)d\mu_{\alpha}(\zeta)\right)^{p_-(B)}\mu_{\alpha}(B^k(z))^{-p(z)}\|\chi_{B^k(z)}\|_{p'(\cdot),w'}^{p(z)}\\
			&\lesssim\left(\frac{1}{\mu_{\alpha}(B^k(z))}\int_{B^k(z)}g(\zeta)d\mu_{\alpha}(\zeta)\right)^{p_-(B)}\mu_{\alpha}(B^k(z))^{p_-(B)-p(z)}\|\chi_{B^k(z)}\|_{p'(\cdot),w'}^{p(z)-p_-(B)}\\
			&\leq\mu_{\alpha}(B^k(z))^{p_-(B)-p(z)}\|\chi_{B^k(z)}\|_{p'(\cdot),w'}^{p(z)-p_-(B)}\frac{1}{\mu_{\alpha}(B^k(z))}\int_{B^k(z)}g(\zeta)^{p_-(B)}d\mu_\alpha(\zeta)\\
			&\lesssim\frac{1}{\mu_{\alpha}(B^k(z))}\int_{B^k(z)}g\chi_{g\geq1}(\zeta)^{p_-(B)}d\mu_\alpha(\zeta)+1\\
			&\lesssim R_k^\alpha g^{p(\cdot)}(z)+1.
			\end{align*}
For the last but one inequality, we also used the following inequality
$$\|\chi_{B^k(z)}\|_{p'(\cdot),w'}^{p(z)-p_-(B)}\leq \max (1, \left \Vert \chi_{\mathbb B}\right \Vert_{p'(\cdot),w'}^{p_+-p_-}).$$			
		\end{proof}
		\begin{lem}\label{cor3}
			Let $p(\cdot)\in\mathcal{P}^{\log}_{\pm}(\mathbb B), k\in(0,\frac{1}{2})$ and $w\in B_{p(\cdot)}$. Then
			\begin{equation*}
			\|R_k^\alpha g\cdot w^{\frac{1}{p(\cdot)}}\|_{p(\cdot)}\lesssim\|g\cdot(R_k^\alpha w)^{\frac{1}{p(\cdot)}}\|_{p(\cdot)}\label{eq23}
			\end{equation*}
			for all non-negative functions $g$ belonging to $L^{p(\cdot)}(R_k^\alpha wd\mu_\alpha).$ 
		\end{lem}
		\begin{proof}
			From Proposition~\ref{lem3}, we have $\sigma=R_k^\alpha w\in A_{p(\cdot)}\subset B_{p(\cdot)}$ because $w\in B_{p(\cdot)}.$ Without loss of generality, we assume that $\|g\|_{p(\cdot),\sigma}=1.$ Thus from Lemma~\ref{rq3} and Lemma~\ref{lem2} we have
			\begin{align*}
			\rho_{p(\cdot)}(R_k^\alpha g\cdot w^{\frac{1}{p(\cdot)}})&=\int_{\mathbb B} (R_k^\alpha g(z))^{p(z)}w(z)d\mu_{\alpha}(z)\\
			&\lesssim\int_{\mathbb B} R_k^\alpha (g^{p(\cdot)})(z)w(z)d\mu_{\alpha}(z)+w(\mathbb B)\\
			&\lesssim\int_{\mathbb B}g(z)^{p(z)} \sigma(z)d\mu_{\alpha}(z)+w(\mathbb B)\\
			&\lesssim1+w(\mathbb B).
			\end{align*}
			For the last inequality, apply Lemma~\ref{le8}.
		\end{proof}
		We still use the notation $\sigma=R_k^\alpha w,$ with $k\in (0, \frac 12).$
		\begin{lem}\label{cor33}
			Let $k\in (0, \frac 12),$  $p(\cdot)\in\mathcal{P}^{\log}_{\pm}(\mathbb B)$ and $w\in B_{p(\cdot)}.$  Then there exists a positive constant $C([w']_{B_{p'(\cdot)}})$ depending  on $[w']_{B_{p'(\cdot)}}$ such that for every non-negative function $g\in L^{p'(\cdot)}(\sigma'),$ we have
			\begin{equation*}
			\|R_k^\alpha g\|_{p'(\cdot),w'}\leq C([w']_{B_{p'(\cdot)}})\|g\|_{p'(\cdot),\sigma'}\label{eq222}
			\end{equation*}
			for all non-negative functions $g\in L^{p'(\cdot)}(\sigma').$ 
		\end{lem}
		\begin{proof}
			Without loss of generality, we assume that  $\|g\|_{p'(\cdot), \sigma'}=1$. From  Proposition \ref{w-w'bis} and Proposition~\ref{lem3}, we have $\sigma'\in A_{p'(\cdot)}\subset B_{p'(\cdot)}$ because $w\in B_{p(\cdot)}$. Thus since $\|g\|_{p'(\cdot),\sigma'}=1$, using Lemma~\ref{rq3} and Lemma~\ref{lem2} we have
			\begin{align*}
			\rho_{p'(\cdot),w'}(R_k^\alpha g)&=\int_{\mathbb B} (R_k^\alpha g(z))^{p'(z)}w'(z)d\mu_{\alpha}(z)\\
			&\lesssim\int_{\mathbb B} R_k^\alpha (g^{p'(\cdot)})(z)w'(z)d\mu_{\alpha}(z)+w'(\mathbb B)\\
			&\lesssim\int_{\mathbb B}g(z)^{p'(z)} R^\alpha_k w'(z)d\mu_{\alpha}(z)+w'(\mathbb B)\\
			&=\int_{\mathbb B}g(z)^{p'(z)}\sigma'(z)\sigma'(z)^{-1} R^\alpha_k w'(z)d\mu_{\alpha}(z)+w'(\mathbb B)\\
			&\leq[w']_{B_{p'(\cdot)}}\int_{\mathbb B}g(z)^{p'(z)}\sigma'(z)d\mu_{\alpha}(z)+w'(\mathbb B)\\
			&\lesssim [w']_{B_{p'(\cdot)}}+w'(\mathbb B).
			\end{align*}
			Indeed, the last inequality follows from Lemma \ref{le8}; for the last but one inequality, using Theorem \ref{equiv-def} and Lemma \ref{le9} for $w'$ in the place of $w,$ we get 
			$$\sigma'(z)^{-1} R^\alpha_k w'(z)=(R_k^\alpha w(z))^{\frac{p'(z)}{p(z)}}R_k^\alpha w'(z)\lesssim [w']_{B_{p'(\cdot)}}$$ 
			because $B^k(z)$ is 'almost' a member of $\mathcal B,$ as it is a subset of the member $B(z, 1-|z|)$ of $\mathcal{B},\,\mu_{\alpha}(B^k(z))\simeq\mu_{\alpha}(B(z, 1-|z|))$ and $w'\in B_{p'(\cdot)}$.
			So by Lemma \ref{le8}, we get  \begin{equation*}
			\|R_k^\alpha g\|_{p'(\cdot),w'}\lesssim C([w']_{B_{p'(\cdot)}}).
			\end{equation*}	 
		\end{proof}
		\begin{lem}\label{lem4}
			Let $p(\cdot)\in\mathcal{P}^{\log}_{\pm}(\mathbb B),\,f\in L^{p(\cdot)}(w),\,k\in \left (0,\frac{1}{2}\right )$ and $w\in B_{p(\cdot)}$. Then
			\begin{equation*}\|R_k^\alpha f (R_k^\alpha w)^{\frac{1}{p(\cdot)}}\|_{p(\cdot)}\leq C([w']_{B_{p'(\cdot)}})\|f\|_{p(\cdot),w}.\label{eq24}\end{equation*}
		\end{lem}
		\begin{proof}
			We still write $\sigma=R^\alpha_k w$. By duality (Proposition \ref{egalm}), there exists a function $g$ satisfying $\|g\|_{p'(\cdot),\sigma'}=1$ and such that
			\begin{equation*}
			\|R_k^\alpha f \|_{p(\cdot),\sigma}\leq 2\int_{\mathbb B}g(z)R^\alpha_kf(z)d\mu_\alpha(z).
			\end{equation*}  
			Next, from Lemma~\ref{lem2}, the H\"older inequality and Lemma~\ref{cor33}, we obtain
			\begin{align*}
			\|R_k^\alpha f \|_{p(\cdot),\sigma}&\lesssim\int_{\mathbb B}f(z)R^\alpha_kg(z)d\mu_\alpha(z)\\
			&\leq 2\|f\|_{p(\cdot),w}\| R^\alpha_k g \|_{p'(\cdot), w'}\\
			&\leq 2C([w']_{B_{p'(\cdot)}})\|f\|_{p(\cdot), w}\|g \|_{p'(\cdot),\sigma'}\\
			&=2C([w']_{B_{p'(\cdot)}})\|f\|_{p(\cdot), w}.
			\end{align*}
			Hence we have the result.
		\end{proof}
		
		\begin{thm}\label{thm2}
			Let $p(\cdot)\in\mathcal{P}^{\log}_{\pm}(\mathbb B)$. If $w\in B_{p(\cdot)},$  there exists a non-negative function $C$ defined on $(0, \infty)$ such that for all $f\in L^{p(\cdot)}(w)$ we have
			\begin{equation*}
			\|m_\alpha f\|_{p(\cdot),w}\leq C([w]_{B_{p(\cdot)}})\|f\|_{p(\cdot),w}.
			\end{equation*}
		\end{thm}
		\begin{proof}
			By Proposition \ref{w-w'}, we have the equality $[w]_{B_{p(\cdot)}}=[w']_{B_{p'(\cdot)}}.$ So			from Lemma~\ref{lem4}, we have
			\begin{equation}\|R_k^\alpha f\cdot(R_k^\alpha w)^{\frac{1}{p(\cdot)}}\|_{p(\cdot)}\leq C([w]_{B_{p(\cdot)}})\|f\|_{p(\cdot),w}.\label{eqt1}\end{equation}
			Hence $R_k^\alpha f\in L^{p(\cdot)}(\sigma).$ Next, since $\sigma \in A_{p(\cdot)}$ by Proposition~\ref{lem3},  Corollary~\ref{rqinv} gives
			\begin{equation}
			\|m_\alpha (R_k^\alpha f)\cdot(R_k^\alpha w)^{\frac{1}{p(\cdot)}}\|_{p(\cdot)}\lesssim\|(R_k^\alpha f)\cdot (R_k^\alpha w)^{\frac{1}{p(\cdot)}}\|_{p(\cdot)}	\label{eqt2}
			\end{equation}
			Hence from \eqref{eqt1} and \eqref{eqt2}, we have $m_\alpha (R_k^\alpha f)\in L^{p(\cdot)}(R_k^\alpha wd\mu_{\alpha})$. 
			Now, Lemma~\ref{cor3}
			gives
			\[\|R_k^\alpha(m_\alpha(R_k^\alpha f))w^{\frac{1}{p(\cdot)}}\|_{p(\cdot)}\lesssim\|m_\alpha (R_k^{\alpha}f)\cdot (R_k^\alpha w)^\frac{1}{p(\cdot)}\|_{p(\cdot)}.\]
			Next, by Proposition \ref{pro1}, there exists a positive constant $C$ such that \[m_\alpha f\leq CR_k^\alpha(m_\alpha(R_k^\alpha f)).\]
			This implies that
			\[\|m_\alpha f\|_{p(\cdot),w}\lesssim\|m_\alpha (R_k^\alpha f)\cdot(R_k^\alpha w)^{\frac{1}{p(\cdot)}}\|_{p(\cdot)}.\]
			Finally, applying \eqref{eqt2} and \eqref{eqt1} successively, we have the result.
		\end{proof}
		\section{A weighted extrapolation Theorem and the proof of the sufficient condition in  Theorem~\ref{th1}}\label{sec-8}
		We are now ready to prove the sufficient condition in Theorem~\ref{th1}, we will adapt the strategy used in \cite{cruz2017extrapolation}.
		\subsection{Preliminary results}
		We recall the $B_1$ class of weights. A weight  $w$ belongs to $B_1$ if 
			\begin{equation}
			[w]_{B_1}:=ess \sup \limits_{z\in \mathbb B} \frac{m_\alpha w(z)}{w(z)}<\infty.
			\end{equation}
In $\mathbb R^n,$ the analogous of the following factorisation theorem was proved for the Muckenhoupt classes $A_p, 1<p<\infty,$ by Jones \cite{jones1980factorization}. 
		\begin{thm}\label{lemextr}
			For a constant exponent $p$ such that $1<p<\infty,$ the following two assertions are equivalent.
			\begin{enumerate}
				\item[1)]
				$w\in B_p;$ 
				\item[2)] there exist $w_1\in B_1$ and $w_2\in B_1$ such that $w=w_1w_2^{1-p}.$
			\end{enumerate}
		\end{thm}
		\begin{proof}
We first show the implication $2)\Rightarrow 1).$ Suppose that $w=w_1w_2^{1-p}$ with $w_1, w_2\in B_1.$ For all $B\in\mathcal{B}$  and $z\in B,$ we have 
			\begin{equation}\frac{1}{\mu_\alpha(B)}\int_B w_id\mu_\alpha\leq[w_i]_{B_1}w_i(z),\quad i=1,2.\label{eqB}\end{equation}
			Thus as $(1-p')(1-p)=1,$ we have $w^{1-p'}=\left(w_1w_2^{1-p}\right)^{1-p'}=w_1^{1-p'}w_2$. So from \eqref{eqB} we have
		\begin{align*}
			    &	\left(\frac{1}{\mu_\alpha(B)}\int_B wd\mu_\alpha\right)\left(\frac{1}{\mu_\alpha(B)}\int_B w^{1-p'}d\mu_\alpha\right)^{p-1} \\
			    &=\left(\frac{1}{\mu_\alpha(B)}\int_B w_1w_2^{1-p}d\mu_\alpha\right)\left(\frac{1}{\mu_\alpha(B)}\int_B w_1^{1-p'}w_2d\mu_\alpha\right)^{p-1}\\
			&\leq[w_1]_{B_1}[w_2]_{B_1}^{p-1}\left(\frac{1}{\mu_\alpha(B)}\int_B w_2d\mu_\alpha\right)^{1-p}\left(\frac{1}{\mu_\alpha(B)}\int_B w_1d\mu_\alpha\right) \\
			&\times\left(\frac{1}{\mu_\alpha(B)}\int_B w_2d\mu_\alpha\right)^{p-1}\left(\frac{1}{\mu_\alpha(B)}\int_B w_1d\mu_\alpha\right)^{-1}\\
			&=[w_1]_{B_1}[w_2]_{B_1}^{p-1}.
	\end{align*}	
\noindent		
Hence $w\in B_p$.
			
We next show the converse implication $2)\Rightarrow 1).$ Suppose that $w\in B_p.$ Set $q=pp'$ and define the operator $S_1$ on the space $\mathcal M$ by
			\[S_1f(z)=w(z)^{\frac{1}{q}}\left(m_\alpha\left(f^{p'}w^{-\frac{1}{p}}\right)(z)\right)^{\frac{1}{p'}}.\]
			By the Minkowski inequality, $S_1$ is sublinear.  Moreover, from the constant exponent version of Theorem \ref{thm2} \cite[Proposition 3]{ref3}, we have
			\[\int_{\mathbb{B}}S_1f(z)^qd\mu_\alpha(z)=\int_{\mathbb{B}}\left(m_\alpha\left(f^{p'}w^{-\frac{1}{p}}\right)(z)\right)^{p}w(z)d\mu_\alpha(z)\lesssim C\left ([w]_{B_p}\right )\int_{\mathbb{B}}f^{q}(z)d\mu_\alpha(z).\]
			In other words, $\|S_1\|_{q}\lesssim \left (C\left ([w]_{B_p}\right )\right )^{\frac{1}{q}}$.\\
 Similarly, denote again $w'=w^{1-p'}\in B_{p'}$ and define the operator $S_2$ on the space $\mathcal M$ by
			\[S_2f(z)=w'(z)^{\frac{1}{q}}\left(m_\alpha\left(f^{p}w'^{-\frac{1}{p'}}\right)(z)\right)^{\frac{1}{p}}.\] 
			By the Minkowski inequality, $S_2$ is also sublinear.  Moreover,
			\[\int_{\mathbb{B}}S_2f(z)^qd\mu_\alpha(z)\lesssim C\left ([w']_{B_{p'}}\right )\int_{\mathbb{B}}f^{q}(z)d\mu_\alpha(z).\]
			In other words, $\|S_2\|_{q}\lesssim C\left (\left ([w']_{B_{p'}}\right )\right )^{\frac{1}{q}}=\left (C\left ([w]_{B_p}\right )\right )^{\frac{1}{q}}$.

		We use the following lemma.
		\begin{lem}\label{S}			
			Set $S=S_1+S_2$ and define the operator $\mathcal{R}$ on $\mathcal M$ by
			\[\mathcal{R}h(z)=\sum_{k=0}^{\infty}\frac{S^kh(z)}{2^k\|S\|_q^k}\]
with $S^0 h=|h|.$ Then 
			\begin{enumerate}
				\item[a)] $|h|\leq \mathcal{R} h$;
				\item[b)] $\|\mathcal{R} h\|_q\leq 2\|h\|_q$;
				\item[c)] $S(\mathcal{R}h)\leq2\|S\|_q\mathcal{R} h$.
			\end{enumerate}
		\end{lem}
		\begin{proof}[Proof of Lemma \ref{S}]
			By the definition of $\mathcal{R}h$ we have $h\leq\mathcal{R}h$. Moreover 
			\begin{align*}
			\|\mathcal{R} h\|_q\leq\sum_{k=0}^{\infty}\frac{\|S^kh\|_q}{2^k\|S\|_q^k}\leq\|h\|_q\sum_{k=0}^{\infty}\frac{1}{2^k}=2\|h\|_q.
			\end{align*}
			Next, the sublinearity of $S$ gives $S(\mathcal{R}h)\leq2\|S\|_q\mathcal{R}h$.
		\end{proof}
			Applying assertion c) of Lemma \ref{S}, we obtain
			\begin{equation}
			w(z)^{\frac{1}{q}}\left(m_\alpha\left((\mathcal{R}h)^{p'}w^{-\frac{1}{p}}\right)(z)\right)^{\frac{1}{p'}}=S_1(\mathcal{R}h)(z)\leq S(\mathcal{R}h)(z)\leq2\|S\|_q\mathcal{R}h(z).\label{eqw1}
			\end{equation}
			Now set $w_2=(\mathcal{R}h)^{p'}w^{-\frac{1}{p}}.$ By \eqref{eqw1}, we have $w_2\in B_1$.   
			
			Similarly, we  have 
			\begin{equation}\label{eqw2}
w'(z)^{\frac{1}{q}}\left(m_\alpha\left(\mathcal{R}h\right )^{p}w'^{-\frac{1}{p'}}(z)\right)^{\frac{1}{p}}=S_2(\mathcal{R}h)(z)\leq S(\mathcal{R}h)\leq2\|S\|_q\mathcal{R}h(z).
\end{equation}
			Now set $w_1=(\mathcal{R}h)^{p}w'^{-\frac{1}{p'}}.$ By \eqref{eqw2}, we have $w_1\in B_1$. Moreover $w_1w_2^{1-p}=w\in B_p.$ This finishes the proof of 
			Lemma~\ref{lemextr}.
			\end{proof}
		
		\begin{lem}\label{lemextr1}
			Let $p(\cdot)\in\mathcal{P}^{\log}_{\pm}(\mathbb B)$ and let $w\in B_{p(\cdot)}.$  We define the operator $R$ on $L^{p(\cdot)}(w)$ by
			\begin{equation*} Rh(x)=\sum_{k=0}^{\infty}\frac{m_\alpha^kh(x)}{2^k\|m_\alpha\|_{L^{p(\cdot)}(w)}^k}\end{equation*} 
			where for $k\geq1, m_\alpha^k={\underbrace { m_\alpha\circ m_\alpha\circ\cdots\circ m_\alpha}_{\textsf{k-times}}}$ and $m_\alpha^0h=|h|$. Then $R$ satisfies the following properties:
			\begin{enumerate}
				\item[a)] $|h|\leq Rh;$
				\item[b)]$R$ is bounded on $L^{p(\cdot)}(w)$ and $\|Rh\|_{p(\cdot),w}\leq2\|h\|_{p(\cdot),w}$;
				\item[c)]$Rh\in B_1$ and $[Rh]_{B_1}\leq2\|m_\alpha\|_{L^{p(\cdot)}(w)}$.
			\end{enumerate}
		\end{lem}
		\begin{proof}
			The proof of assertions a) and b) are the same as for assertions a) and b) of Lemma \ref{S}. Here, we use the sublinearity of $m_\alpha.$ 
			\\Finally, by the definition of $Rh,$  we have
			\begin{align*}
			m_\alpha (Rh)(x)&\leq\sum_{k=0}^{\infty}\frac{m_\alpha^{k+1}h(x)}{2^k\|m_\alpha\|_{L^{p(\cdot)}(w)}^k}\\
			&\leq2\|m_\alpha\|_{L^{p(.)}(w)}\sum_{k=0}^{\infty}\frac{m_\alpha^{k+1}h(x)}{2^{k+1}\|m_\alpha\|_{L^{p(.)}(w)}^{k+1}}\\
			&\leq2\|m_\alpha\|_{L^{p(.)}(w)}Rh(x).
			\end{align*}
			 Thus $Rh\in B_1$ and $[Rh]_{B_1}\leq2\|m_\alpha\|_{L^{p(.)}(w)}$.	
		\end{proof}
		\begin{lem}\label{lemextr2}
			Let $p(\cdot)\in\mathcal{P}^{\log}_{\pm}(\mathbb B)$ and $w\in B_{p(\cdot)}$. Define the operator $H$ on $L^{p'(\cdot)}$ by
			\[Hh=\mathcal{R}'\left(hw^{\frac{1}{p(\cdot)}}\right)w^{-\frac{1}{p(\cdot)}}\]where 
			\[\mathcal{R}'g(x)=\sum_{k=0}^{\infty}\frac{m_\alpha^kg(x)}{2^k\|m_\alpha\|_{L^{p'(.)}(w')}^k}.\]Then
			\begin{enumerate}
				\item[a)] $|h|\leq Hh;$
				\item[b)]$H$ is bounded on $L^{p'(\cdot)}$ and $\|Hh\|_{p'(\cdot)}\leq2\|h\|_{p'(\cdot)}$;
				\item[c)]$Hh\cdot w^{\frac{1}{p(\cdot)}}\in B_1$ and $[Hh\cdot w^{\frac{1}{p(\cdot)}}]_{B_1}\leq2\|m_\alpha\|_{L^{p'(.)}(w')}$.
			\end{enumerate}
		\end{lem}
		\begin{proof}
			The proof is the same as for Lemma \ref{lemextr1}. We replace $p(\cdot)$ by $p'(\cdot)$ and $w\in B_{p(\cdot)}$ by $w'\in B_{p'(\cdot)}.$ The property $p'(\cdot)\in \mathcal{P}^{\log}_{\pm}(\mathbb B)$ comes from Remark \ref{p-p'}.
		\end{proof}
		\subsection{A weighted extrapolation theorem}
		We denote by $\mathcal{F}$  a family of couples of non-negative measurable functions.	We are now ready to state and prove the following weighted variable extrapolation theorem.	
		\begin{thm}\label{thextrax}
			Suppose that for some constant exponent $p_0>1,$ there exists a  function $C: (0, \infty)\to(0, \infty)$  such that for all $v\in B_{p_0}$ and  $(F,G)\in\mathcal{F}$, we have
			\begin{equation}
			\int_{\mathbb B}F(x)^{p_0}v(x)d\mu_\alpha(x)\leq C([v]_{B_{p_0}})\int_{\mathbb B}G(x)^{p_0}v(x)d\mu_\alpha(x).\label{eqext}
			\end{equation} 
			Then given $p(\cdot)\in\mathcal{P}^{\log}_{\pm}(\mathbb B)$ and $w\in B_{p(\cdot)}$,  we have
			\begin{equation}\|F\|_{p(\cdot),w}\leq 16\times 4^{-\frac 1{p_0}}\left (C([v]_{B_{p_0}})\right )^{\frac 1{p_0}}\|G\|_{p(\cdot),w}\end{equation}
			for all $(F,G)\in\mathcal{F}$ and $F\in L^{p(\cdot)}(w)$
		\end{thm}
		\begin{proof}
			We use the  technique of Cruz-Uribe in \cite[Theorem~$2.6$]{cruz2017extrapolation}.
			Let $(F,G)\in\mathcal{F}$. If $\|F\|_{p(\cdot),w}=0$ we have the result. Otherwise, $\|F\|_{p(\cdot),w}>0$ and hence $\|G\|_{p(\cdot),w}>0,$ because if $\|G\|_{p(\cdot),w}=0$, then $G=0$ a.e. and by  \eqref{eqext} we will have $F=0$ a.e. Henceforth, we assume $0<\|F\|_{p(\cdot),w}<\infty$ and $0<\|G\|_{p(\cdot),w}<\infty.$ Define 
			\[h_1=\frac{F}{\|F\|_{p(\cdot),w}}+\frac{G}{\|G\|_{p(\cdot),w}},\]
			then $\|h_1\|_{p(\cdot),w}\leq2$ and so $h_1\in L^{p(\cdot)}(w)$.
			
			Since $F\in L^{p(\cdot)}(w)$, by duality (Proposition \ref{egalm}), there exists $h_2\in L^{p'(\cdot)}$ such that $\|h_2\|_{p'(\cdot)}=1$ and 
			\begin{equation}
			\|F\|_{p(\cdot),w}\leq 2\int_{\mathbb B}Fw^{\frac{1}{p(\cdot)}}h_2d\mu_{\alpha}\leq 2\int_{\mathbb B}F(Hh_2)w^{\frac{1}{p(\cdot)}}d\mu_{\alpha}\label{eqex1}
			\end{equation}
			where the latter inequality comes from assertion $a)$ of Lemma~\ref{lemextr2}.
			
			Set $\gamma =\frac 1{p'_0}.$ By the usual H\"older inequality, we have
			\begin{align}
			\int_{\mathbb B}F(Hh_2)w^{\frac{1}{p(\cdot)}}d\mu_{\alpha}&=\int_{\mathbb B}F(Rh_1)^{-\gamma}(Rh_1)^{\gamma}(Hh_2)w^{\frac{1}{p(\cdot)}}d\mu_{\alpha}\nonumber\\
			&\leq I_1^{\frac{1}{p_0}}I_2^{\frac{1}{p'_0}},\label{eqex2}
			\end{align}
			where
			$$I_1:=\int_{\mathbb B}F^{p_0}(Rh_1)^{1-p_0}(Hh_2)w^{\frac{1}{p(\cdot)}}d\mu_{\alpha}$$
			and
			$$I_2:=\int_{\mathbb B}(Rh_1)(Hh_2)w^{\frac{1}{p(\cdot)}}d\mu_{\alpha}.$$
			In addition, from Lemma~\ref{lemextr1} and Lemma~\ref{lemextr2} respectively,  $R$ is bounded on  $L^{p(\cdot)}(w)$ and $H$ is bounded on $L^{p'(\cdot)}.$ Thus by the H\"older inequality, assertions $b)$ of Lemma~\ref{lemextr1} and Lemma~\ref{lemextr2}, we have
			\[I_2\leq2\|Rh_1\|_{p(\cdot),w}\|Hh_2\|_{p'(\cdot)}\leq 8 \|h_1\|_{p(\cdot),w}\|h_2\|_{p'(\cdot)} \leq 16.\] 
			By the definition of $h_1$ and assertion a) of Lemma~\ref{lemextr1}, we have $$\frac{\varphi}{\|\varphi\|_{p(\cdot),w}}\leq h_1\leq Rh_1 $$ 
			for $\varphi\in \{F, \hskip 1truemm G\}.$ Next, by the H\"older inequality and assertion b) of Lemma~\ref{lemextr2}, we have
			\begin{align*}
			I_1
			&\leq \int_{\mathbb B} F^{p_0}(\zeta)\left(\frac{F(\zeta)}{\|F\|_{p(\cdot),w}}\right)^{1-p_0}H(\zeta)h_2(\zeta)w^{\frac{1}{p(\zeta)}}d\mu_{\alpha}(\zeta)\\
			&=\|F\|_{p(\cdot),w}^{p_0-1}\int_{\mathbb B} F(\zeta)H(\zeta)h_2(\zeta)w^{\frac{1}{p(\zeta)}}d\mu_{\alpha}(\zeta)\\
			&\leq2\|F\|_{p(\cdot),w}^{p_0-1}\|F\|_{p(\cdot),w}\|Hh_2\|_{p'(\cdot)}
			\\&\leq4\|F\|_{p(\cdot),w}^{p_0}
			\\&<\infty.
			\end{align*} 
			Since $Rh_1\in B_1$ and $(Hh_2)w^{\frac{1}{p(\cdot)}}\in B_1$ by Lemma~\ref{lemextr1} and Lemma~\ref{lemextr2} respectively, it follows from Theorem~\ref{lemextr} that $v:=(Rh_1)^{1-p_0}\left(Hh_2w^{\frac{1}{p(\cdot)}}\right)\in B_{p_0}.$
			Hence by \eqref{eqext} and the same argument as above, we have
			\begin{align*}
			I_1&=\int_{\mathbb B} F^{p_0}(Rh_1)^{1-p_0}(Hh_2)w^{\frac{1}{p(\cdot)}}d\mu_{\alpha}\\
			&\leq C([v]_{B_{p_0}})\int_{\mathbb B} G^{p_0}(Rh_1)^{1-p_0}(Hh_2)w^{\frac{1}{p(\cdot)}}d\mu_{\alpha}\\
			&\leq C([v]_{B_{p_0}})\int_{\mathbb B} G^{p_0}\left(\frac{G}{\|G\|_{p(\cdot),w}}\right)^{1-p_0}(Hh_2)w^{\frac{1}{p(\cdot)}}d\mu_{\alpha}\\
			&=C([v]_{B_{p_0}})\|G\|_{p(\cdot),w}^{p_0-1}\int_{\mathbb B} G(Hh_2)w^{\frac{1}{p(\cdot)}}d\mu_{\alpha}\\
			&\leq2C([v]_{B_{p_0}})\|G\|_{p(\cdot),w}^{p_0-1}\|G\|_{p(\cdot),w}\|Hh_2\|_{p'(\cdot)}
			\\&\leq4C([v]_{B_{p_0}})\|G\|_{p(\cdot),w}^{p_0}.
			\end{align*}
			Thus from \eqref{eqex1} and \eqref{eqex2}, we have the result.
		\end{proof}
		\subsection{The end of the proof of the sufficient condition in Theorem \ref{th1}}
We prove the following proposition.
		\begin{prop}\label{prop1}
			Let $p(\cdot)\in\mathcal{P}^{\log}_{\pm}(\mathbb B)$ and $w\in B_{p(\cdot)}.$ Then $P_\alpha^+$ is a continuous operator on $L^{p(\cdot)}(w).$ Consequently, the Bergman projector $P_\alpha$ extends to a continuous operator on $L^{p(\cdot)}(w).$
		\end{prop}
		\begin{proof}
			We call again $\mathcal{C}_c(\mathbb B)$ the space   of continuous functions of compact support in $\mathbb B$ and we take $\mathcal{F}=\{(P^+_\alpha f, \vert f\vert):\,\;f\in \mathcal{C}_c(\mathbb B)\}$. We recall from Proposition \ref{dense} that $\mathcal{C}_c(\mathbb B)$ is a dense subspace in $L^{p(\cdot)} (w).$
			
			Let $p_0$ be an arbitrary constant exponent greater than $1.$  Let $v\in B_{p_0}.$ By Theorem~\ref{Bek_class}, for every $f\in \mathcal{C}_c(\mathbb B),$  we have  
			\[\int_{\mathbb B}(P^+_\alpha f)^{p_0}vd\mu_{\alpha}\leq C\left ([v]_{B_{p_0}}\right )\int_{\mathbb B}|f|^{p_0}vd\mu_{\alpha}.\]
Thus by Theorem~\ref{thextrax}, for all $f\in \mathcal{C}_c(\mathbb B),$ we have
			\[\|P^+_\alpha f\|_{p(\cdot),w}\leq 16\times 4^{-\frac 1{p_0}}\left (C([v]_{B_{p_0}})\right )^{\frac 1{p_0}}\|f\|_{p(\cdot),w}.\]
			We conclude by density.
		\end{proof}


\end{document}